\documentclass[12pt]{amsart}

\usepackage{amsmath,amssymb,amsthm}
\usepackage{amscd,times,fullpage}
\newtheorem{thm}{Theorem}
\newtheorem{prop}{Proposition}
\newtheorem{lem}{Lemma}
\newtheorem{cor}{Corollary}

\theoremstyle{remark}
\newtheorem{rem}{Remark}

\theoremstyle{definition}
\newtheorem{defn}{Definition}
\newtheorem{ex}{Example}

\newcommand{\C}{\mathbb{ C}}

\newcommand{\N}{\mathbb{ N}}

\newcommand{\Q}{\mathbb{ Q}}
\newcommand{\R}{\mathbb{ R}}

\newcommand{\rk}{\operatorname{rk}}

\newcommand{\Sym}{\operatorname{Sym}}
\newcommand{\Sh}{\operatorname{Sh_\lambda}}
\newcommand{\Z}{\mathbb{ Z}}

\newcommand{\sg}{\mathrm{sign}}
\newcommand{\rw}{\mathrm{w}}

\DeclareMathOperator{\dummygg}{\mathfrak{g}}
\renewcommand{\gg}{\dummygg}
\DeclareMathOperator{\II}{\mathfrak{I}}
\DeclareMathOperator{\hh}{\mathfrak{h}}
\DeclareMathOperator{\TT}{\mathfrak{t}}

\DeclareMathOperator{\LL}{\mathfrak{L}}
\DeclareMathOperator{\ad}{ad}
\DeclareMathOperator{\Ad}{Ad}

\title{Equivariant complex structures on homogeneous
spaces and their cobordism classes}
\author{Victor M.~Buchstaber}
\address{Steklov Mathematical Institute, Russian Academy
of Sciences, Gubkina Street 8, 119991 Moscow, Russia}
\email{buchstab@mi.ras.ru}
\address{School of Mathematics, University of Manchester, Oxford Road,
Manchester M13~9PL, England}
\email{Victor.Buchstaber@manchester.ac.uk}
\author{Svjetlana Terzi\'c}
\address{Faculty of Science, University of Montenegro,
Cetinjski put bb, 81000 Podgorica, Montenegro}
\email{sterzic@cg.ac.yu}
\date{\today ; MSC 2000: primary 57R77, 22F30; secondary 57R20, 14M15 }

\begin{document}

\maketitle

\begin{abstract}
We consider compact homogeneous spaces $G/H$, where $G$ is a compact
connected Lie group and $H$ is its closed connected subgroup of
maximal rank. The aim of this paper is to provide an effective
computation of the universal toric genus for the complex, almost
complex and stable complex structures which are invariant under the
canonical left action of the maximal torus $T^k$ on $G/H$. As it is
known, on $G/H$ we may have many such structures and the
computations of their toric genus in terms of fixed points for the
same torus action give the constraints on possible collections of
weights for the corresponding representations of $T^k$ in the
tangent spaces at the fixed points, as well as on the signs at these
points.  In that context, the effectiveness is also approached due
to an explicit description of the relations between the weights and
signs for an arbitrary couple of such structures. Special attention
is devoted to the structures which are invariant under the canonical
action of the group $G$. Using classical results, we obtain an
explicit description of the weights and signs in this case. We
consequently obtain an expression for the cobordism classes of such
structures in terms of coefficients of the formal group law in
cobordisms, as well as in terms of Chern numbers in cohomology.
These computations require no information on the cohomology ring of
the manifold $G/H$, but, on their own, give important relations in
this ring. As an application we provide an explicit formula for the
cobordism classes and characteristic numbers of the flag manifolds
$U(n)/T^n$, Grassmann manifolds $G_{n,k}=U(n)/(U(k)\times U(n-k))$
and some particular interesting examples.

\end{abstract}

\section{introduction}

The homogeneous spaces $G/H$, where $G$ is a compact connected Lie
group and $H$ is its closed connected subgroup of maximal rank are
classical manifolds, which play significant role in  many areas of
mathematics. They can be characterized as homogeneous spaces with
positive Euler characteristic. Our interest in these manifolds is
related to the well known  problem in cobordism theory to find
representatives  in cobordism classes that have reach geometric
structure and carry many non-equivalent stable complex structures.
Let us mention~\cite{B} that all compact homogeneous spaces $G/H$
where $H$ is centralizer of an element of odd order admit an
invariant almost complex structure. Furthermore, if $H$ is a
centralizer of a maximal torus in $G$ then $G/H$ admits an invariant complex
structure which gives rise to a unique invariant Kaehler-Einstein
metric and moreover all homogeneous complex spaces are algebraic.
Besides that by~\cite{WG}, the stationary subgroup $H$ of any
homogeneous complex space $G/H$ can be realized  as the fixed point
set under the action of some finite group of inner automorphisms of
$G$ and vice versa. This shows that these spaces can be made off
generalized symmetric spaces  what results in the  existence of
finite order symmetries on them. Our interest in the study of the
homogeneous spaces $G/H$ with positive Euler characteristic is also
stimulated by the well known relations between the cohomology rings of
these spaces and the deep problems in the theory of representations
and combinatorics (see, for example~\cite{F}). These problems are
formulated in terms of different additive basis in the cohomology
rings for $G/H$ and multiplicative rules related to that basis. We
hope that the relations in the cohomology rings $H^{*}(G/H,\Z )$
obtained from the calculation of the universal toric genus of the
manifold $G/H$ may lead to new results in that direction.

In the paper~\cite{Novikov-67}, which opened a new stage in the
development of the cobordism theory, S.~P.~Novikov proposed a method
for the description of the fixed points for actions of groups on
manifolds, based on the formal group law for geometric cobordisms.
That paper  rapidly stimulated active research work which brought
significant results. These results in the case of $S^1$-actions are
mainly contained in the
papers~\cite{Gusein-Zade-71-Izv},~\cite{Gusein-Zade-71-MZ},~\cite{Krichever-74},~\cite{Krichever-76}
and also in~\cite{Krichever-90}. Our approach to this problem uses
the results on the universal toric genus, which was introduced
in~\cite{BR} and described in details in~\cite{IMRN}. Let us note
that the formula for the universal toric genus in terms of fixed
points is a generalization of Krichever's
formula~\cite{Krichever-74} to the case of stable complex manifolds.
For the description of the cobordism classes of manifolds in terms
of their Chern numbers we appeal to the Chern-Dold character theory,
which is developed in~\cite{Buchstaber}.

The universal toric genus can be constructed for any even
dimensional manifold $M^{2n}$ with a given torus action and stable
complex structure which is equivariant under the torus action.
Moreover, if the set of isolated fixed points for this action is
finite than the universal toric genus can be localized, which means
that it can be explicitly written through the weights and the signs
at the fixed points for the representations that gives arise from
the given torus action.

The construction of the toric genus is reduced to the computation of
Gysin homomorphism of $1$ in complex cobordisms for fibration whose
fiber is $M^{2n}$ and the base is classifying space of the torus.
The problem of the localization of Gysin homomorphism is very known
and it was studied by many authors, starting with 60-es of the last
century. In~\cite{BR} and~\cite{IMRN} is obtained explicit answer
for this problem in the terms of the torus action on tangent bundle
for $M^{2n}$. The history of this problem is presented also in these
papers.

If consider compact homogeneous space $G/H$ with $\rk G = \rk H = k$,
we have on it the canonical action $\theta$ of the maximal torus
$T^{k}$ for $H$ and $G$, and any $G$-invariant almost complex structure
on $G/H$ is compatible with this action. Besides that, all fixed
points for the action $\theta$ are isolated, so one can apply
localization formula to compute universal toric genus for this
action and any invariant almost complex structure. Since, in this
case, we consider almost complex structures, all fixed points in the
localization formula  are going to have sign $+1$. We prove that the
weights for the action $\theta$ at different fixed points can be
obtained  by the action of the Weyl group $W_{G}$ up to the action of
the Weyl group $W_{H}$ on the weights for $\theta$ at identity fixed
point. On this way we get an  explicit formula for the cobordism
classes of such spaces in terms of the weights at the fixed point
$eH$. This formula also shows that the cobordism class for $G/H$
related to an invariant almost complex structure can be computed
without information about cohomology for $G/H$.

We obtain also the explicit formulas, in terms of the weights at
identity fixed point, for the cohomology characteristic numbers for
homogeneous spaces of positive Euler characteristic endowed with an
invariant almost complex structure. We use further that the
cohomology characteristic numbers $s_\omega, \; \omega=(i_1, \ldots,
i_n)$, and classical Chern numbers $c^\omega=c_1^{i_1}\cdots
c_n^{i_n}$ are related by some standard relations from the theory of
symmetric polynomials. This fact together with the obtained formulas
for the characteristic numbers $s_\omega(\tau(G/H))$ proves that the
classical Chern numbers $c^\omega(\tau(G/H))$ for the almost complex homogeneous
spaces  can be computed without information on
their cohomology. It also gives an explicit way for the computation
of the classical Chern numbers.

In studying invariant almost complex structures on compact
homogeneous spaces $G/H$ of positive Euler characteristic we appeal
to the theory developed in~\cite{BH} which describes such
structures in terms of complementary roots for $G$ related to $H$. In
that context by an invariant almost complex structure on $G/H$ we
assume the structure that is invariant under the canonical action
of $G$ on $G/H$.

In this paper we go further with an application of the universal toric
genus and consider the almost complex structures on $G/H$ for which
we require to be invariant only under the canonical action $\theta$
of the maximal torus $T$ on $G/H$, as well as the stable complex
structures equivariant under this torus action. We prove generally that for an arbitrary
$M^{2n}$ the weights for any two such structures  differ
only by the signs. We also prove that the sign difference in the
weights determines, up to common factor $\pm 1$, the difference of
the signs at the  fixed points related to any two such structures.

We provide an application of our results by obtaining explicit
formula for the cobordism class and top cohomology characteristic
number of the flag manifolds $U(n)/T^n$ and Grassmann manifolds
$G_{n,k}=U(n)/(U(k)\times U(n-k))$ related to the standard complex
structures. We want to emphasize that, our method when applying to
the flag manifolds and Grassmann manifolds gives the description of
their cobordism classes and characteristic numbers using the
technique of {\it divided difference operators}. Our method also
makes possible to compare cobordism classes that correspond to the
different invariant almost complex structures on the same
homogeneous space. We illustrate that on  the space
$U(4)/(U(1)\times U(1)\times U(2))$, which is firstly given
in~\cite{BH} as an example of homogeneous space  that admits two
different invariant complex structures.

We also compute the universal toric genus of the  sphere
$S^{6}=G_{2}/SU(3)$ related to unique $G_2$-invariant almost complex
structure for which is known not to be complex~\cite{BH}. We prove
more, that this structure is also unique almost complex structure
invariant under the canonical action $\theta$ of the maximal torus
$T^{2}$.

This paper comes out from the first part of our work where we mainly
considered invariant almost complex structures on homogeneous spaces
of positive Euler characteristic. It has continuation which is going
to deal with the same questions, but related to the  stable complex
structures equivariant under given torus action on homogeneous
spaces of positive Euler characteristic.

The authors are grateful to A.~Baker and A.~Gaifullin for useful
comments.

\section{Universal toric genus}
We will recall the results from~\cite{BR},~\cite{BPR} and~\cite{IMRN}.

\subsection{General setting.}In general setting one considers
$2n$-dimensional manifold $M^{2n}$ with a given smooth action
$\theta$ of the torus $T^{k}$. We say that $(M^{2n}, \theta,
c_{\tau})$ is {\it equivariant stable complex} if it admits
$\theta$-equivariant stable complex structure $c_{\tau}$. This means
that there exists $l\in \N$ and complex vector bundle $\xi$ such that
\begin{equation}\label{scdp}
c_{\tau}\colon \tau (M^{2n})\oplus \R ^{2(l-n)}\longrightarrow \xi
\end{equation}
is real isomorphism and the composition
\begin{equation}\label{scd}
r(t)\colon\xi\stackrel{c_\tau^{-1}}{\longrightarrow}\tau(M^{2n})\oplus
\mathbb{R}^{2(l-n)}\stackrel{d\theta(t)\oplus I}{\longrightarrow}
\tau(M^{2n})\oplus\mathbb{R}^{2(l-n)}\stackrel{c_\tau}{\longrightarrow}\xi
\end{equation}
is a complex transformation for any $t\in T^{k}$.

If there exists $\xi$ such that $c_{\tau}\colon \tau
(M^{2n})\longrightarrow \xi$ is an isomorphism, i.~e.~ $l=n$, then
$(M^{2n},\theta ,c_{\tau})$ is called {\it almost complex}
$T^k$-manifold.

Denote by $\Omega _{U}^{*} [[u_1,\ldots, u_{k}]]$ an algebra of
formal power series over $\Omega _{U}^{*} = U^{*}(pt)$. It is well
known~\cite{Novikov} that $U^{*}(pt) = \Omega _{U}^{*} = \Z
[y_{1},\ldots ,y_{n},\ldots ]$, where $\dim y_{n}=-2n$. Moreover, as
the generators for $\Omega _{U}^{*}$ over the rationales, or in
other words for $\Omega _{U}^{*}\otimes \Q$, can be taken the family
of cobordism classes $[\C P^{n}]$ of the complex projective spaces.

When given a $\theta$-equivariant stable complex structure
$c_{\tau}$ on $M^{2n}$,  we can always choose $\theta$-equivariant
embedding $i\colon M^{2n}\to \R ^{2(n+m)}$, where $m>n$, such that
$c_{\tau}$ determines, up to natural equivalence, a
$\theta$-equivariant complex structure $c_{\nu}$ on the normal
bundle $\nu(i)$ of $i$. Therefore, one can define the universal
toric genus for $(M^{2n},\theta, c_{\tau})$ in complex cobordisms,
see~~\cite{BR},~\cite{IMRN}.

We want to note that, in the case when $c_{\tau}$ is almost complex
structure, a universal toric genus for $(M^{2n},\theta ,c_{\tau})$
is completely defined in terms of the action $\theta$ on tangent
bundle $\tau (M^{2n})$.

The universal toric genus for $(M^{2n}, \theta, c_{\tau})$ could be
looked at as an element in algebra $\Omega _{U}^{*} [[u_1,\ldots,
u_{k}]]$. It is defined with
\begin{equation}\label{uto}
\Phi (M^{2n}, \theta, c_{\tau}) =  [M^{2n}] + \sum _{|\omega|>0}
[G_{\omega}(M^{2n})]u^{\omega} \; ,
\end{equation}
where $\omega = (i_1,\ldots ,i_k)$ and $u^{\omega} =
u_{1}^{i_1}\cdots u_{k}^{i_k}$.\\
Here by $[M^{2n}]$ is denoted the complex cobordism class of the
manifold $M^{2n}$  with stable complex structure $c_{\tau}$, by
$G_{\omega}(M^{2n})$ is denoted the stable complex manifold obtained
as the total space of the fibration $G_{\omega} \to B_{\omega}$ with
fiber $M$. The base $B_{\omega} = \prod _{j=1}^{k}B_{j}^{i_j}$,
where $B_{j}^{i_j}$ is Bott tower, i.~e.~$i_j$-fold iterated
two-sphere bundle over $B_{0}=pt$. The base $B_{\omega}$ satisfies
$[B_{\omega}]=0$, $|\omega|>0$, where $|\omega| =
\sum\limits_{j=1}^{k} i_{j}$.

For the universal toric genus of homogeneous space of positive Euler characteristic we prove the following.

\begin{lem}\label{L1}
Let $M^{2n}=G/H$, where $G$ is compact connected Lie group and $H$
its closed connected subgroup of maximal rank. Denote by $\theta$
the canonical action of the maximal torus $T^k$ on $M^{2n}$ and let
$c_{\tau}$ be $G$-equivariant stable complex structure on $M^{2n}$.
Then the universal toric genus $\Phi (M^{2n},\theta, c_{\tau})$
belongs to the image of homomorphism $Bj^{*}\colon U^{*}(BG) \to
U^{*}(BT^{k})$ which is induced by the embedding $T^{k}\subset G$.
\end{lem}
\begin{proof}
According to its construction, the universal toric genus $\Phi
(M^{2n},\theta ,c_{\tau})$ is equal to $p_!(1)$ and belongs to
$U^{-2n}(BT^{k})$, where
\[
p : ET^{k}\times _{T^k} M^{2n}\to BT^{k} \; .
\]
Using that the action $\theta$ is induced by the left action of the
group $G$ and looking at the commutative diagram
\begin{eqnarray}
ET^{k} & \rightarrow & EG\nonumber\\
\downarrow & & \downarrow\nonumber\\
BT^{k} & \rightarrow & BG\nonumber
\end{eqnarray}
we obtain the proof due to the fact that Gysin homomorphism is
functorial for the bundles that can be connected with commutative
diagram.
\end{proof}

\subsection{The action  with isolated fixed points.}

We first introduce, following~\cite{BR}, the general notion of the
{\it sign} at isolated fixed point. Let we are on $M^{2n}$ given an
equivariant stable complex structure $c_{\tau}$.

We assume $\R ^{2(l-n)}$ in~\eqref{scdp} to be endowed with
canonical orientation. Under this assumption the real
isomorphism~\eqref{scdp} defines an orientation on $\tau(M^{2n})$,
since $\xi$ is canonically oriented by an existing complex
structure.

On the other hand, if $p$ is the isolated fixed point, the
representation $r_{p}\colon T^k\to GL(l,\C)$ associated
to~\eqref{scd} produces the decomposition  of the fiber $\xi
_{p}\cong \C ^{l}$ as $\xi _{p}\cong \C ^{l-n}\oplus \C ^{n}$. In
this decomposition $r_{p}$ acts trivially on $\C ^{l-n}$ and without
trivial summands on $\C ^{n}$. Note that $\C^{n}$ inherits here the
complex structure from $\xi _{p}$ which defines an orientation of
$\C^{n}$. This together leads to the following definition.

\begin{defn}
 The  $\sg (p)$ at isolated fixed point $p$ is
$+1$ if the map
\[
\tau_p(M^{2n})\stackrel{I\oplus\hspace{.1ex}0}{\longrightarrow}
\tau_p(M^{2n})\oplus \mathbb{R}^{2(l-n)} \stackrel{c_{\tau,p}}
{\longrightarrow}\xi_p \cong\mathbb{C}^n \oplus\mathbb{C}^{l-n}
\stackrel{\pi}{\longrightarrow} \mathbb{C}^n \; ,
\]
preserves orientation. Otherwise, $\sg (p)$ is $-1$.
\end{defn}

\begin{rem}\label{cmx}
Note that for an almost complex $T^{k}$-manifold $M^{2n}$, it
directly follows from the definition that  $\sg (p)= +1$ for any
isolated fixed point.
\end{rem}

If an action $\theta$ of $T^{k}$ on $M^{2n}$ has only isolated fixed
points, then it is proved that toric genus for $M^{2n}$  can be
completely described using just local data at the fixed
points,~\cite{BR},~\cite{IMRN}.

Namely, let $p$ again be an isolated fixed point. Then the non
trivial summand of $r_p$ from~\eqref{scd} gives rise to the
tangential representation of $T^{k}$ in $GL(n, \C)$. This
representation decomposes into $n$ non-trivial one-dimensional
representations of $r_{p,1}\oplus \ldots \oplus r_{p,n}$ of $T^{k}$.
Each of the representations $r_{p,j}$ can be written as
\[
r_{p,j}(e^{2\pi i x_{1}}, \ldots ,e^{2\pi i x_{k}})v= e^{2\pi i
\langle\Lambda_{j}(p), {\bf x}\rangle }v \; ,
\]
for some $\Lambda_{j}(p) = (\Lambda^{1}_{j}(p),\ldots
,\Lambda^{k}_{j} (p)) \in \Z ^{k}$, where ${\bf x}=(x_{1},\ldots
,x_{k})\in \R ^{k}$ and $\langle \Lambda_{j}(p), {\bf x}\rangle =
\sum\limits _{l=1}^{k}\Lambda^{l}_{j}(p)x_{l}$. The sequence $\{
\Lambda_{1}(p),\ldots ,\Lambda_{n}(p)\}$ is called {\it the weight
vector} for representation $r_{p}$ in the fixed point $p$.

\begin{rem}
Since $p$ is an isolated fixed point none of the couples of  weights $\Lambda_{i}(p)$, $1\leq i\leq n$, have common integer factor.
\end{rem}

\begin{thm}\label{weights2}
Let $c_{\tau}^{1}$ and $c_{\tau}^{2}$ be the stable complex
structures on $M^{2n}$ equivariant under the given action $\theta$
of the torus $T^{k}$ on the manifold $M^{2n}$ with isolated fixed
points.
\begin{itemize}
\item The weights $\Lambda_{i}^{1}(p)$ and $\Lambda_{i}^{2}(p)$,
$1\leqslant i\leqslant n$, of an action $\theta$ at the fixed point
$p$ corresponding to the structures $c_{\tau}^{1}$ and $c_{\tau}^{2}$
are related with
\begin{equation}\label{swgeneral}
\Lambda_{i}^{2}(p) = a_{i}(p)\Lambda_{i}^{1}(p), \;  \mbox{where} \;
a_{i}(p)=\pm 1 \; .
\end{equation}
\item The signs $\sg(p)_{1}$ and $\sg(p)_{2}$ at the fixed point $p$
corresponding to the structures $c_{\tau}^{1}$ and $c_{\tau}^{2}$ are
related with
\begin{equation}\label{signsgeneral}
\sg(p)_{2} = \epsilon \cdot \prod_{i=1}^{n}a_{i}(p)\cdot \sg(p)_{1} \; ,
\end{equation}
where $\epsilon = \pm 1$ depending if $c_{\tau}^{1}$ and
$c_{\tau}^{2}$ define of $\tau(M^{2n})$ the same orientation or not.
Here $a_{i}(p)$ are such that
$\Lambda^{2}_{i}(p)=a_{i}(p)\Lambda^{1}_{i}(p)$, for the weights
$\Lambda^{1}_{i}(p)$ and $\Lambda^{2}_{i}(p)$ of an action $\theta$
related to the structures $c_{\tau}^{1}$ and $c_{\tau}^{2}$.
\end{itemize}
\end{thm}
\begin{proof}
Let $p$ be an isolated fixed point of an action $\theta$ of the
torus $T^{k}$ on the manifold $M^{2n}$. If it is  on  $M^{2n}$ given
the $\theta$-equivariant stable complex structure $c_{\tau_{1}}$
than in the neighborhood of $p$, the tangential representation of
$T^{k}$ in $GL(n,\C )$ assigned to an action $\theta$ and structure
$c_{\tau_{1}}$ decomposes into the sum of non-trivial
one-dimensional representations $r_{p,1}\oplus\ldots \oplus
r_{p,n}$. Any other stable complex structure $c_{\tau_{2}}$ which is
equivariant under the given action $\theta$ commutes with  each of
the one-dimensional representations $r_{p,i}$, $1\leqslant
i\leqslant n$. Therefore, the one-dimensional summands in which
decomposes the tangential representation of $T^{k}$, assigned to an
action $\theta$ and structure $c_{\tau_{2}}$, are $r_{p,i}$ or it's
conjugate $\overline{r_{p,i}}$, $1\leqslant i\leqslant n$. This
implies that the relations between the weights for an action
$\theta$ related to the two different stable $\theta$-equivariant
stable complex structures are given by the
formula~\eqref{swgeneral}.

To prove the second statement of the Theorem let us note that the
sign at the fixed point $p$ of some $\theta$-equivariant 
stable complex structure $c_{\tau}$  is  determined by the
orientations of the real two-dimensional subspaces in which decomposes summand $\C^{n}$ of
$\xi _{p}= \C^{n}\oplus \C^{l-n}$. That decomposition is obtained using the decomposition of the 
tangential representation of $T^{k}$  determined by the action $\theta$ and the structure
$c_{\tau}$. Therefore, by~\eqref{swgeneral} it follows that the
relation between the signs at the fixed point  for the given torus
action related to the two equivariant stable complex structures is
given by~\eqref{signsgeneral}.
\end{proof}

\begin{rem}
We want to point that Theorem~\ref{weights2} gives that, under
assumption that  manifold $M^{2n}$ admits $\theta$-equivariant
stable complex structure, the signs at the fixed points for any
other $\theta$-equivariant stable complex structure  are completely  determined by an orientation
that structure defines on $M^{2n}$ and by t weights at fixed points. In
other words, when passing from the existing to some other
$\theta$-equivariant stable complex structure, the "correction" of
the sign at arbitrary fixed point is completely determined by the
"correction" of the weights at that fixed point up to some common
factor $\epsilon=\pm 1$ which points on the difference of
orientations on $M^{2n}$ that these two structures define.
\end{rem}

\subsection{Formal group low.} Let $F(u, v)= u + v + \sum
\alpha _{ij}u^{i}v^{j}$ be {\it the formal group for complex
cobordism} \cite{Novikov-67}. The corresponding power system $\{
[w](u)\in \Omega ^{*}[[u]] : w\in \Z \}$ is uniquely defined with
$[0](u)=0$ and $[w](u) = F(u, [w-1])(u)$, for $w\in \Z$. For ${\bf
w}=(w_1,\ldots ,w_{k})\in {\Z}^{k}$ and ${\bf u}=(u_1,\ldots, u_k)$
one defines $\bf{[w](u)}$ inductively with ${\bf [w](u)} =[w](u)$
for $k=1$ and
\[
{\bf[w](u)} = F_{q=1}^{k}[w_{q}](u_{q})=F(F_{q=1}^{k-1}[w_{q}](u_{q}),
[w_{k}](u_{k})) \; ,
\]
for $k\geqslant 2$. Then for toric genus of the action $\theta$
with isolated fixed points the following localization formula holds,
which is first formulated in~\cite{BR} and proved in details
in~\cite{IMRN}.

\begin{thm}\label{EGFP}
If the action $\theta$ has a finite set $P$  of isolated fixed
points then
\begin{equation}\label{tw}
\Phi (M^{2n}, \theta, c_{\tau}) = \sum _{p\in P} \sg (p) \prod
_{j=1}^{n}\frac{1}{[\Lambda_{j}(p)]({\bf u})} \; .
\end{equation}
\end{thm}

\begin{rem}
Theorem~\ref{EGFP} together with  formula~\eqref{uto} gives that
\begin{equation}\label{order}
\sum _{p\in P} \sg (p) \prod
_{j=1}^{n}\frac{1}{[\Lambda_{j}(p)]({\bf u})}= [M^{2n}] + \LL ({\bf u})\; ,
\end{equation}
where $\LL ({\bf u}) \in \Omega _{U}^{*} [[u_{1},\ldots ,u_{k}]]$
and $\LL ({\bf 0}) = 0$. In this way Theorem~\ref{EGFP} gives that
all summands in the left hand side of~\eqref{order} have order $n$
in $0$.
\end{rem}

\begin{rem}
As we will make it explicit further, the fact that after making the
sum, all singularities in formula~\eqref{tw} should disappear, gives
constraints on the weights and signs at the fixed points. Note also,
that formula~\eqref{tw} gives an expression for the cobordism class
$[M^{2n}]$ in terms of the weights and signs at fixed points.
\end{rem}

\subsection{Chern-Dold character.}\label{CDS}  We show further how one,
together with Theorem~\ref{EGFP}, can use the notion of Chern-Dold
character in cobordisms in order to obtain an expression  for
cobordism class $[M^{2n}]$ in terms of the characteristic numbers for $M^{2n}$,
as well as the relations on the weights and signs at fixed points.
In review of the basic definitions and results on Chern character we
follow~\cite{Buchstaber}.

Let $U^{*}$ be the theory of unitary cobordisms.

\begin{defn}
The Chern-Dold character  for a topological space $X$ in the theory
of unitary cobordisms $U^{*}$ is a ring homomorphism
\begin{equation}
ch _{U} : U^{*}(X) \to H^{*}(X, \Omega _{U}^{*}\otimes \Q) \ .
\end{equation}
\end{defn}

Recall that the Chern-Dold character as a multiplicative
transformation of cohomology theories is uniquely defined by the
condition that for $X=(pt)$ it gives canonical inclusion $\Omega
_{U}^{*}\to \Omega _{U}^{*}\otimes \Q$.

The Chern-Dold character splits into composition
\begin{equation}\label{CD}
ch _{U} : U^{*}(X)\to H^{*}(X, \Omega _{U}^{*}(\Z ))\to H^{*}(X,
\Omega _{U}^{*}\otimes \Q ) \; .
\end{equation}
The ring $\Omega _{U}^{*}(\Z )$ in \eqref{CD} is firstly described
in~\cite{Buchstaber}. It is a subring of $\Omega _{U}^{*}\otimes \Q$
generated by the elements from $\Omega _{U}^{-2n}\otimes \Q$ having
integer Chern numbers. It is equal to
\[
\Omega _{U}^{*}(\Z )=\Z [b_1,\ldots , b_n,\ldots ] \; ,
\]
where $b_n = \frac{1}{n+1}[\C P^{n}]$.

The Chern character leaves  $[M^{2n}]$ invariant, i.~e.~
\[
ch_{U}([M^{2n}])= [M^{2n}]\; ,
\]
and $ch_{U}$ is the homomorphism of $\Omega _{U}^*$-modules

It follows from the its description~\cite{Buchstaber} that the
Chern-Dold character $ch _{U} : U^{*}(X)\to H^{*}(X, \Omega
_{U}^{*}(\Z ))$ as a multiplicative transformation of the cohomology
theories is given by the series
\[
ch _{U}u = h(x)=\frac{x}{f(x)},\quad \mbox{where}\quad f(x)= 1+\sum
_{i=1}^{\infty} a_{i}x^{i} \quad \mbox{and}\quad a_{i}\in \Omega
_{U}^{-2i}(\Z) \; .
\]
Here $u=c_1^U(\eta) \in U^2(\mathbb{C}P^\infty)$ and $x=c_1^H(\eta)
\in H^2(\mathbb{C}P^\infty,\mathbb{Z})$ denote the first Chern classes
of the universal complex line bundle $\eta \to \mathbb{C}P^\infty$.

From the construction of Chern-Dold character it follows also the
equality
\begin{equation}\label{charact}
ch_{U}[M^{2n}] = [M^{2n}] = \sum _{\|\omega \|=n} s_{\omega}(\tau
(M^{2n}))a^{\omega}\; ,
\end{equation}
where $\omega = (i_{1},\ldots ,i_{n}), \; \| \omega \| = \sum
_{l=1}^{n} l\cdot i_{l}$ and $a ^{\omega} = a_{1}^{i_1}\cdots
a_{n}^{i_{n}}$. Here the numbers $s_{\omega}(\tau
(M^{2n})),{\|\omega \|=n}$ are {\it the cohomology characteristic
numbers} of $M^{2n}$ and they correspond to the  cohomology tangent
characteristic classes of $M^{2n}$.

If on $M^{2n}$ is given torus action $\theta$ of $T^{k}$ and stable
complex structure $c_{\tau}$ which is $\theta$-equivariant, then the
Chern character of its toric genus is
\begin{equation}\label{cheq}
ch_{U}\Phi (M^{2n},\theta, c_{\tau}) = [M^{2n}] + \sum _{| \xi |
>0}[G_{\xi}(M^{2n})](ch_{U}{\bf u})^{\xi}\, ,
\end{equation}
where $ch_{U}{\bf u} = (ch_{U}u_1, \ldots, ch_{U}u_k),\;\; ch_{U}u_i
= \frac{x_i}{f(x_i)}$\, and\, $\xi = (i_1, \ldots,i_k),\; | \xi | =
i_1+ \cdots+i_k$.

We have that $F(u, v)=g^{-1}\left( g(u)+g(v)\right)$, where $g(u)=u
+ \sum\limits_{n>0}\frac{1}{n+1}[\C P^{n}]u^{n+1}$ (see
\cite{Novikov-67}) is {\it the logarithm of the formal group} $F(u,
v)$ and $g^{-1}(u)$ is the exponent of $F(u, v)$, that is the
function inverse to the series $g(u)$. Using that $ch_{U}g(u)=
g(ch_{U}(u))= g\big( \frac{x}{f(x)} \big) =x$ (see
\cite{Buchstaber}), we obtain $g^{-1}(x)=\frac{x}{f(x)}$ and
$ch_{U}F(u_1, u_2)=\frac{x_1+x_2}{f(x_1+x_2)}$ and therefore
\[ ch_{U}[\Lambda _{j}(p)](u) = \frac{\langle \Lambda _{j}(p),
{\bf x}\rangle}{f(\langle \Lambda _{j}(p), {\bf
x}\rangle)}\; .  \]
Applying these results to Theorem~\ref{EGFP} we get
\begin{equation}\label{chs}
ch_{U}\Phi(M^{2n},\theta, c_{\tau })= \sum _{p\in P}\sg (p)\prod
_{j=1}^{n}\frac{f(\langle \Lambda _{j}(p), {\bf x}\rangle)}{\langle
\Lambda _{j}(p),{\bf x}\rangle} \; .
\end{equation}
From~\eqref{cheq} and~\eqref{chs} it follows  that
\begin{equation}\label{cc}
 \sum _{p\in P}\sg (p)\prod _{j=1}^{n}\frac{f(\langle \Lambda _{j}(p),
{\bf x}\rangle)}{\langle \Lambda _{j}(p), {\bf x}\rangle}= [M^{2n}]
+ \sum _{| \xi | >0}[G_{\xi}(M^{2n})](ch_{U}{\bf u})^{\xi} \, .
\end{equation}

\begin{ex}\label{n=2}
Let us take $M^{2}=\C P^{1}=U(2)/\big(U(1)\times U(1)\big)$. We have
the action $\theta$ of $T^2$ on $\C P^{1}$ with two fixed points. The
weights related to the standard complex structure $c_\tau$, are
$(x_1-x_2)$ and $(x_2-x_1)$. By equation~\eqref{chs} we obtain
that the Chern character of the universal toric genus for
$(\mathbb{C}P^1,\theta,c_\tau)$ is given by the series
\[ ch_{U}\Phi(\mathbb{C}P^1,\theta,c_\tau) = \frac{f(x_1-x_2)}{x_1-x_2}
+ \frac{f(x_2-x_1)}{x_2-x_1} = 2\sum_{k=0}^{\infty} a_{2k+1}
(x_1-x_2)^{2k} \, . \] By  equation~\eqref{cheq} we obtain
\[ [\C P^1] + \sum_{i+j>0}[G_{i,j}(\C P^1)]\frac{x_1^ix_2^j}{f(x_1)^if(x_2)^j}
= 2a_1 + 2\sum\limits_{k=1}^{\infty} a_{2k+1} (x_1-x_2)^{2k} \, . \]
Thus, $[\C P^1]=2a_1$ and $\sum\limits_{i+j=n}[G_{i,j}(\C P^1)]=0$
for any $n>0$. Moreover,
\begin{align*}
[G_{i,j}(\C P^1)]\; & \sim \;0, \; \text{if }\;i+j=2k+1,\\
[G_{i,2k-j}(\C P^1)]\; & \sim \;(-1)^i2\binom{2k}{i}a_{2k+1}, \;
k>0,
\end{align*}
where ``$\sim$'' is equality to the elements decomposable in
$\Omega_U(\mathbb{Z})$.

Note that the subgroup $S^{1}=\{(t_1,t_2)\in T^2 | t_2=1\}$ acts also on $\C P^{1}$ with the same fixed points
as for the action of $T^2$, but the weights are going to be $x$ and $-x$. This action is, as an example,  given in~\cite{BR}.
\end{ex}

If in the left hand side of the equation~\eqref{cc} we put $t{\bf x}$ instead
of ${\bf x}$ and then multiply it with $t^{n}$ we obtain the
following result.

\begin{prop}\label{cobc}
The coefficient for $t^{n}$ in the series in $t$
\[
\sum _{p\in P}\sg (p)\prod _{j=1}^{n}\frac{f(t\langle \Lambda
_{j}(p), {\bf x}\rangle)}{\langle \Lambda _{j}(p), {\bf x}\rangle}
\]
represents the complex cobordism class $[M^{2n}]$.
\end{prop}

\begin{prop}\label{coeffzero}
The coefficient for $t^{l}$ in the series in $t$
\[
\sum _{p\in P}\sg (p)\prod _{j=1}^{n}\frac{f(t\langle \Lambda
_{j}(p), {\bf x}\rangle)}{\langle \Lambda _{j}(p), {\bf x}\rangle}
\]
is equal to zero for $0\leqslant l\leqslant n-1$.
\end{prop}

\section{Torus action on homogeneous spaces with positive Euler
characteristic. }

Let $G/H$ be a compact homogeneous space of positive Euler
characteristic. It means that $G$ is a compact connected Lie group
and $H$ its connected closed subgroup, such that $\rk G = \rk H$.
Let $T$ be the maximal common torus for $G$ and $H$. There is
canonical action $\theta$ of $T$ on $G/H$ given by $t(gH)=(tg)H$,
where $t\in T$ and $gH\in G/H$. Denote by $N_{G}(T)$ the normalizer
of the torus $T$ in $G$. Then $W_{G} = N_{G}(T)/T$ is the Weyl group
for $G$. For the set  of fixed points  for the action $\theta$
we prove the following.

\begin{prop}\label{fixed}
The set of fixed points under the canonical action $\theta$ of $T$ on
$G/H$ is given by $(N_{G}(T))~\cdot~H$.
\end{prop}
\begin{proof}
It is easy to see  that $gH$ is fixed point for $\theta$ for any
$g\in N_{G}(T)$. If $gH$ is the fixed point under the canonical
action of $T$ on $G/H$ then $t(gH) = gH$ for all $t\in T$. It
follows that $g^{-1}tg \in H$ for all $t\in T$, i.~e.~ $g^{-1}Tg
\subset H$. This gives that $g^{-1}Tg$ is a maximal torus in $H$
and, since any two maximal toruses in $H$ are conjugate, it follows
that $g^{-1}Tg = h^{-1}Th$ for some $h\in H$. Thus,
$(gh)^{-1}T(gh)=T$ what means that $gh\in N_{G}(T)$. But,
$(gh)H=gH$, what proves the statement.
\end{proof}

Since $T \subset N_{G}(T)$ leaves $H$ fixed, the following Lemma is
direct implication of the  Proposition~\ref{fixed}.

\begin{lem}\label{wfixed}
The set of fixed points under the canonical action $\theta$ of $T$ on
$G/H$ is given by $W_{G}\cdot H$.
\end{lem}

Regarding the number of fixed points, it holds the following.

\begin{lem}\label{number}
The number of fixed points under the canonical action $\theta$ of $T$
on $G/H$ is equal to the Euler characteristic $\chi (G/H)$.
\end{lem}
\begin{proof}
Let $g, g^{'}\in N_{G}(T)$ are representatives of the same fixed
point. Then $g^{'}g^{-1}\in H$ and $g^{-1}Tg = T =
(g^{'})^{-1}Tg^{'}$, what gives that $g^{'}g^{-1}Tg(g^{'})^{-1} = T$
and, thus,  $g^{'}g^{-1}\in N_{H}(T)$. This implies that the number
of fixed points is equal to
\[
\Big{\|} \frac{N_{G}(T)}{N_{H}(T)} \Big{\|} = \frac{\| \frac{N_{G}
(T)} {T}\|}{\| \frac{N_{H}(T)}{T} \|} = \frac{\| W_{G}\|}{\| W_{H}
\|}  = \chi (G/H) \; .
\]
The last equality is classical result related to equal ranks
homogeneous spaces, see~\cite{Onishchik}.
\end{proof}

\begin{rem}
The proof of Lemma~\ref{number} gives that the set of fixed
points under the canonical action $\theta$  of $T$ on $G/H$ can be
obtained as an orbit of $eH$ by the action of the Weyl  group
$W_{G}$ up to the action of the Weyl group $W_{H}$.
\end{rem}

\section{The weights at the fixed points.}

Denote  by $\gg$, $\hh$ and $\TT$ the Lie algebras for $G$, $H$ respectively and
$T=T^k$, where $k=\rk G=\rk H$. Let $\alpha
_{1},\ldots, \alpha _{m}$ be the roots for $\gg$ related to $\TT$,
where $\dim G=2m+k$.  Recall that the roots for $\gg$ related to
$\TT$ are the weights for the adjoint representation $\Ad _{T}$  of
$T$ which is given  with $\Ad _{T}(t) = d_{e}\ad (t)$, where $\ad
(t)$ are inner automorphisms of $G$ defined by the elements $t\in T$.
One can always choose the roots for $G$  such that $\alpha
_{n+1},\ldots, \alpha _{m}$ gives the roots for $\hh$ related to
$\TT$, where $\dim H=2(m-n)+k$. The roots $\alpha _{1},\ldots,
\alpha  _{n}$ are called the {\it complementary} roots for $\gg$
related to $\hh$. Using root decomposition for $\gg$ and $\hh$ it
follows that $T_{e}(G/H) \cong \gg _{\alpha _{1}}^{\C}\oplus \ldots
\oplus \gg _{\alpha _{n}}^{\C}$, where by $\gg _{\alpha _{i}}$ is
denoted the root subspace defined with the root $\alpha _{i}$ and
$T_{e}(G/H)$ is the tangent space for $G/H$ at the $e\cdot H$.
 It is obvious
that $\dim _{\R} G/H=2n$.
\subsection{Description of the invariant almost complex structures.}
Assume we are given  an invariant almost complex structure $J$ on
$G/H$. This means that $J$ is invariant under the canonical action
of $G$ on $G/H$. Then according to the paper~\cite{BH}, we can say
the following.
\begin{itemize}\label{ac}
\item Since $J$ is invariant it commutes with adjoint representation
$Ad_{T}$ of the torus T. This implies that $J$ induces the complex
structure on each complementary root subspace $\gg _{\alpha
_{1}},\ldots ,\gg _{\alpha _{n}}$. Therefore, $J$ can be completely
described by the root system $\varepsilon _{1} \alpha _{1},\ldots
,\varepsilon _{n}\alpha _{n}$, where we take $\varepsilon _{i}= \pm
1$ depending  if $J$ and adjoint representation $Ad_{T}$ define the
same orientation on $\gg _{\alpha _{i}}$ or not,  $1\leqslant
i\leqslant n$. The roots $\varepsilon _{k}\alpha _{k}$ are called
{\it the roots of the almost complex structure} $J$.
\item If we assume $J$ to be integrable, it follows that it can be
chosen the ordering on the canonical coordinates of $\TT$ such that
the roots  $\varepsilon _{1} \alpha _{1},\ldots ,\varepsilon _{n}
\alpha _{n}$ which define $J$ make the closed system of positive
roots.
\end{itemize}
Let us assume that $G/H$ admits an invariant almost complex
structure. Consider the  isotropy representation $I_{e}$ of $H$ in
$T_{e}(G/H)$ and let it decomposes into $s$ {\it real irreducible
representations} $I_e = I_e^1 + \ldots +I_e^s$. Then it is proved
in~\cite{BH} that  $G/H$ admits exactly $2^s$ invariant almost
complex structures. Because of completeness we recall the proof of
this fact shortly here. Consider the decomposition of $T_{e}(G/H)$
\[
T_{e}(G/H) = \II _{1}\oplus \ldots \oplus\II _{s}
\]
such that the restriction of $I_e$ on $\II _{i}$ is $I_{e}^i$. The
subspaces $\II _1,\ldots ,\II_{s}$ are  invariant under $T$ and
therefore each of them is the sum of some root subspaces, i.e. $\II
_{i} = \gg _{\alpha _{i_1}}\oplus\ldots \oplus \gg _{\alpha
_{i_j}}$, for some complementary roots $\alpha _{i_1},\ldots ,\alpha
_{i_j}$. Any linear transformation that commutes with $I_e$ leaves
each of $\II _i$ invariant. Since, by assumption $G/H$ admits
invariant almost complex structure, we have at least one linear
transformation without real eigenvalue that commutes with $I_e$.
This implies that the commuting field for each of $I_{e}^i$ is the
field of complex numbers and, thus, on each $\II _{i}$ we have
exactly two invariant complex structures.

\begin{rem}
Note that this consideration shows that the numbers $\varepsilon
_1,\ldots ,\varepsilon _n$ that define an invariant almost complex
structure may not vary independently.
\end{rem}

\begin{rem}
In this paper we consider almost complex structures on $G/H$ that
are invariant under the canonical action of the group $G$, what, as
we remarked, imposes some relations on $\varepsilon _1,\ldots
,\varepsilon _n$. If we do not require $G$-invariance, but just
$T$-invariance, we will have more degrees of freedom on $\varepsilon
_1,\ldots ,\varepsilon _n$. This paper is going to have continuation,
where, among the other, the case of $T$-invariant structures will be
studied.
\end{rem}

\begin{ex}
Since the isotropy representation for $\C P^{n}$ is irreducible over
the reals, it follows that on $\C P^{n}$ we have only two invariant
almost complex structures, which are actually the standard complex
structure and its conjugate.
\end{ex}

\begin{ex}
The flag manifold $U(n)/T^{n}$ admits $2^{m}$ invariant almost
complex structures, where $m=\frac{n(n-1)}{2}$. By~\cite{BH} only two
of them, conjugate to each other, are integrable.
\end{ex}

\begin{ex}
As we already mentioned, the $10$-dimensional manifold
$M^{10}=U(4)/(U(1)\times U(1)\times U(2))$ is the first example of
homogeneous space, where we have an existence of two non-equivalent
invariant complex structures, see~\cite{BH}. We will, in the last
section of this paper, also describe cobordism class of $M^{10}$ for
these structures.
\end{ex}

\subsection{The weights at the fixed points.}
We fix now an invariant almost complex structure $J$ on $G/H$ and we
want to describe the weights of the canonical action $\theta$ of $T$
on $G/H$ at the fixed points of this action. If $gH$ is the fixed
point for the action $\theta$, then we have a linear map
$d_{g}\theta (t)\colon T_{g}(G/H)\to T_{g}(G/H)$ for all $t\in T$.
Therefore, this action gives rise to the complex representation
$d_{g}\theta$ of $T$ in $(T_{g}(G/H), J)$.

The weights for this representation at identity fixed point  are
described in~\cite{BH}.

\begin{lem}
The weights for the representation $d_{e}\theta$ of $T$ in
$(T_{e}(G/H),J)$ are given by the roots of  an invariant almost
complex structure $J$.
\end{lem}
\begin{proof}
Let us, because of clearness, recall the proof. The inner
automorphism $\ad (t)$, for $t\in T$ induces the map $\overline{\ad
}(t) : G/H\to G/H$ given with $\overline{\ad }(t)(gH) = t(gH)t^{-1} =
(tg)H$. Therefore, $\theta (t) =  \overline{\ad} (t)$ and, thus,
$d_{e}\theta (t) = d_{e}\overline{\ad } (t)$ for any $t\in T$. This
directly gives that the weights  for $d_{e}\theta$ in $(T_{e}(G/H),
J)$ are the roots that define $J$.
\end{proof}

For an arbitrary fixed point we prove the following.

\begin{thm}\label{weights}
Let $gH$ be the fixed point for the canonical action $\theta$ of $T$
on $G/H$. The weights of the induced representation $d_{g}\theta $
of $T$ in $(T_{g}(G/H), J)$ can be obtained from the weights of the
representation $d_{e}\theta$ of $T$ in $(T_{e}(G/H), J)$ by the
action of the Weyl group $W_{G}$ up to the action of the Weyl group
$W_{H}$.
\end{thm}
\begin{proof}
Note  that Lemma~\ref{wfixed} gives that an arbitrary fixed point
can be written as $\rw H$ for some $\rw \in W_{G}/W_{H}$. Fix $\rw
\in W_{G}/W_{H}$ and denote by $l(\rw )$ the action of $\rw$ on
$G/H$, given by $l(\rw)gH=(\rw g)H$ and by $\ad (\rw )$ the inner
automorphism of $G$ given by $\rw$.

We observe that $\theta \circ \ad (\rw )  = \ad (\rw )\circ \theta
$ and therefore $d_{e}\theta \circ d_{e}\ad (\rw ) = d_{e}\ad (\rw )\circ
d_{e}\theta$. This implies that  the weights for $d_{e}\theta \circ
d_{e}\ad (\rw )$ we get by the action of $d_{e}\ad (\rw)$ on the
weights for $d_{e}\theta$. From the other side $\theta (\ad (\rw
)t)gH = ({\rw}^{-1}t\rw g)H = (l({\rw}^{-1})\circ \theta (t)\circ
l(\rw ))gH$ what implies that $d_{e}(\theta \circ \ad (\rw )) =
d_{\rw}l({\rw}^{-1})\circ d_{\rw}\theta \circ d_{e}l(\rw )$. This
gives that if, using the map $d_{\rw}l_{{\rw}^{-1}}$, we lift the
weights for $d_{\rw}\theta$ from $T_{\rw}(G/H)$ to $T_{e}(G/H)$, we
get that they coincide with the weights for $d_{e}\theta \circ
d_{e}\ad (\rw)$. Therefore, the weights for $d_{\rw}\theta$ we can
get by the action of the element $\rw$ on the weights for
$d_{e}\theta$.
\end{proof}

\section{The cobordism classes of homogeneous spaces with positive
Euler characteristic}

\begin{thm}
Let $G/H$ be a homogeneous space of compact connected Lie group such
that $\rk G=\rk H = k$, $\dim G/H=2n$ and consider the canonical
action $\theta$ of maximal torus $T=T^k$ for $G$ and $H$  on $G/H$.
Assume we are given an invariant almost complex structure $J$ on
$G/H$. Let  $\Lambda_{j}=\varepsilon _{j} \alpha _{j}$, $1\leqslant j\leqslant
n$, where $\varepsilon _{1} \alpha _{1},\ldots ,\varepsilon
_{n}\alpha _{n}$ are the complementary roots of $G$ related to $H$
which define an invariant almost complex structure $J$. Then the
toric genus for $(G/H, J)$ is given with
\begin{equation}
\Phi (G/H, J) = \sum _{\rw \in W_{G}/W_{H}} \prod _{j=1}^{n}\frac{1}
{ [\rw (\Lambda_{j})]({\bf u})} \ .
\end{equation}
\end{thm}
\begin{proof}
Rewriting   Theorem~\ref{EGFP}, since all fixed points have sign
$+1$, we get that the toric genus for $(G/H, J)$  is
\begin{equation}\label{EQH}
\Phi (G/H, J) = \sum _{p\in P} \prod _{j=1}^{n}\frac{1}{[\Lambda
_{j}(p)]({\bf u})} \; ,
\end{equation}
where $P$ is the set of isolated fixed points and $\{\Lambda
_{1}(p),\ldots ,\Lambda _{n}(p)\}$ is the weight vector of the
representation for $T$  in  $T_{p}(G/H)$  associated to an action
$\theta$. By Theorem~\ref{number}, the set of fixed points $P$
coincides with the orbit of the action of $W_{G}/W_{H}$ on $eH$ and
also by Theorem~\ref{weights} the set of weight vectors at fixed
points coincides with the orbit of the action of $W_{G}/W_{H}$ on
the weight vector $\Lambda$ at $eH$.  The result follows if we put
this data into formula~\eqref{EQH}.
\end{proof}

\begin{cor}\label{ch-tor-hom}
The Chern-Dold character of the toric genus for homogeneous space
$(G/H, J)$ is given with
\begin{equation}\label{CC}
ch_{U}\Phi(G/H, J) = \sum _{\rw \in W_{G}/W_{H}}\prod _{j=1}^{n}
\frac{f (\langle \rw (\Lambda_{j}), {\bf x}\rangle )}{\langle \rw
(\Lambda_{j}),{\bf x}\rangle} \; ,
\end{equation}
where $f(t) = 1+ \sum\limits_{i\geqslant 1}a_{i}t^{i}$ for $a_{i}
\in \Omega _{U}^{-2i}(\Z), \; {\bf x}=(x_1,\ldots ,x_{k})$ and by
$\langle \Lambda _{j}, {\bf x} \rangle=\sum\limits_{l=1}^k
\Lambda_j^lx_l$ is denoted the weight vector $\Lambda _{j}$ of
$T^k$-representation at $e\cdot H$.
\end{cor}

\begin{cor}\label{chom}
The cobordism class for $(G/H, J)$ is given as the coefficient for
$t^{n}$ in the series in $t$
\begin{equation}
\sum _{\rw \in W_{G}/W_{H}}\prod _{j=1}^{n} \frac{f(t\langle \rw
(\Lambda_{j}), {\bf x}\rangle )}{\langle \rw (\Lambda_{j}) ,{\bf
x}\rangle} \; .
\end{equation}
\end{cor}

\begin{rem}
Since the weights of different invariant almost complex structures
on the fixed homogeneous space $G/H$ differ only by sign, 
Corollary~\ref{chom} provides the way for comparing cobordism
classes of two such structures on $G/H$ without having their cobordism
classes explicitly computed.
\end{rem}

\section{Characteristic numbers of homogeneous spaces
with positive Euler characteristic.}

\subsection{Generally about stable complex manifolds.} Let $M^{2n}$ be
an equivariant stable complex manifold whose given action $\theta$ of
the torus $T^{k}$ on $M^{2n}$ has only isolated fixed points. Denote
by $P$ the set of fixed points for $\theta$ and set $t_{j}(p) =
\langle \Lambda _{j}(p), {\bf x}\rangle$, where $\{ \Lambda _{j}(p),
\; j=1,\ldots,n\}$ are the weight vectors of the representation of
$T^{k}$ at a fixed point $p$ given by the action $\theta$ and ${\bf
x}=(x_{1},\ldots ,x_{k})$.

Set
\begin{equation}\label{fdecomp}
\prod_{i=1}^{n}f(t_i)=1+\sum f_\omega(t_1, \ldots, t_n)a^\omega \; .
\end{equation}
Using this notation  Proposition~\ref{coeffzero} could be
formulated in the following way.

\begin{prop}\label{coeffzero1}
 For any $\omega$ with  $0\leqslant \|\omega\| \leqslant (n-1)$ we have that
\[
\sum _{p\in P}\sg (p)\cdot \frac{f_\omega(t_1(p), \ldots,
t_n(p))}{t_1(p) \cdots t_n(p)} = 0 \; .
\]
\end{prop}

Note that Proposition~\ref{coeffzero1} gives the strong
constraints on the set of signs $\{\sg (p)\}$ and the set of weights
$\{\Lambda _{j}(p)\}$ at fixed points for a manifold with a given
torus action and related equivariant stable complex structure.  For
example for $\omega = (i_1,\ldots ,i_n)$ such that $i_{k}=1$ for
exactly one $k$ such that $1\leqslant k\leqslant n-1$ and $i_{j}=0$
for $j\neq k$ it gives that the signs and the weights at fixed points have to satisfy
the following relations.

\begin{cor}\label{second}
\[
\sum _{p\in P}\sg (p)\cdot \frac{\sum\limits _{i=1}^{n}t_{j}^{k}(p)}{t_1(p)
\cdots t_n(p)} = 0 \; ,
\]
where $0\leqslant k\leqslant n-1$.
\end{cor}

As we already mentioned in~\eqref{charact} the cobordism class for
$M^{2n}$ can be represented as
\[
[M^{2n}]= \sum_{\| \omega \|= n}s_{\omega}(\tau (M^{2n}))a^{\omega} \; ,
\]
where $\omega = (i_{1},\ldots ,i_{n})$, $\| \omega \| = \sum _{l=1}^{n}
l\cdot i_{l}$ and $a ^{\omega} = a_{1}^{i_1}\cdots a_{n}^{i_{n}}$.

If the given action $\theta$ of $T^{k}$ on $M^{2n}$ is with isolated
fixed points, the coefficients $s_{\omega}(\tau (M^{2n}))$ can be
explicitly described using Proposition~\ref{cobc} and
expression~\eqref{fdecomp}.

\begin{thm}\label{sthm}
Let $M^{2n}$ be an equivariant stable complex manifold whose given
action $\theta$ of the $T^{k}$ has  only isolated fixed points.
Denote by $P$  the set of fixed points for $\theta$ and set
$t_{j}(p) = \langle \Lambda _{j}(p), {\bf x}\rangle$, where $\Lambda
_{j}(p)$ are the weight vectors of the representation of $T^{k}$ at fixed points
given by the action $\theta$ and ${\bf x}=(x_{1},\ldots x_{k})$.
Then for $\|\omega\|=n$
\begin{equation}\label{somega}
s_{\omega}(\tau (M^{2n}))=\sum _{p\in P}\sg (p)\cdot
\frac{f_\omega(t_1(p), \ldots, t_n(p))}{t_1(p) \cdots t_n(p)} \; .
\end{equation}
\end{thm}

\begin{ex}\label{ex2}
\[ s_{(n,0,\ldots,0)}(\tau (M^{2n}))=
\sum _{p\in P}\sg (p)\; . \]
\end{ex}

\begin{ex}\label{ex3}
\[ s_{(0,\ldots,0,1)}(\tau (M^{2n}))=s_n(M^{2n})=
\sum _{p\in P}\sg (p) \frac{\sum\limits_{j=1}^n t_j^n(p)}{t_1(p)
\cdots t_n(p)} \; . \]
\end{ex}

\begin{rem}\label{int}
Note that the left hand side of~\eqref{somega} in 
Theorem~\ref{sthm}  is an  integer number $s_{\omega}(\tau
(M^{2n}))$ while the right hand side is a rational function in
variables $x_1, \ldots, x_k$. So this theorem imposes strong
restrictions on the sets of signs $\{\sg (p)\}$ and weight vectors
$\{\Lambda_j (p)\}$ at the fixed points.
\end{rem}

\subsubsection{On existence of more stable complex structures.}
Let assume that manifold $M^{2n}$ endowed with torus action $\theta$
with isolated fixed points admits $\theta$-equivariant stable
complex structure $c_{\tau}$. Denote by $\Lambda_{1}(p),\ldots
,\Lambda _{n}(p)$ the weights for an action
$\theta$ at fixed points $p\in P$ and by $\sg(p)$ the signs at the
fixed points related to $c_{\tau}$. Let further $t_{j}(p) =  \langle
\Lambda_{j}(p), {\bf x}\rangle$, $1\leqslant j\leqslant n$. Then
Theorem~\ref{weights2}, Proposition~\ref{coeffzero1} and
Theorem~\ref{sthm} give the following necessary condition for the
existence of  another  $\theta$-equivariant stable complex
structure on $M^{2n}$..

\begin{prop}\label{neccond}
If $M^{2n}$ admits an other $\theta$-equivariant stable complex
structure $(M^{2n},c^{'}_{\tau},\theta)$ then  there exist
$a_{i}(p)=\pm 1$, where $p\in P$ and $1\leqslant i\leqslant n$, such
that the following conditions are satisfied:
\begin{itemize}
\item for any $\omega$ with  $0\leqslant \|\omega\| \leqslant (n-1)$
\begin{equation}
\sum _{p\in P}\sg (p)\frac{f_\omega(a_{1}(p)t_1(p), \ldots,
a_{n}(p)t_n(p))}{t_1(p) \cdots t_n(p)} = 0 \; .
\end{equation}
\item for any $\|\omega\|=n$
\begin{equation}\label{somegaad}
\sum _{p\in P}\sg (p)
\frac{f_\omega(a_{1}(p)t_1(p), \ldots, a_{n}(p)t_n(p))}{t_1(p) \cdots t_n(p)}
\end{equation}
is an integer number.
\end{itemize}
\end{prop}

As a special case we get analogue of Corollary~\ref{second}.

\begin{cor}
If $M^{2n}$ admits an other $\theta$-equivariant stable complex
structure $(M^{2n},c^{'}_{\tau},\theta)$ then  there exist
$a_{j}(p)=\pm 1$, where $p\in P$ and $1\leqslant j\leqslant n$, such
that
\begin{equation}\label{secondstable}
\sum _{p\in P}\sg (p)\frac{\sum\limits _{i=1}^{n}(a_{j}(p))^{k}t_{j}^{k}(p)}{t_1(p)
\cdots t_n(p)} = 0 \; ,
\end{equation}
for $1\leqslant k\leqslant n-1$.
\end{cor}

\begin{rem}
Note that the relations~\eqref{secondstable} only  for $k$ being odd
give the constraints on the existence of the second stable complex
structure.   In the same spirit it follows from Example~\ref{ex3}
that if $n$ is even than the characteristic number $s_{n}(M^{2n})$
are the same  for all stable complex structures on $M^{2n}$
equivariant under the fixed torus action $\theta$. For $n$ being
odd, as  Example~\ref{cps} and Subsection~\ref{M10}  will show,
these numbers may be different.
\end{rem}

\subsection{Homogeneous spaces of positive Euler characteristic
and with invariant almost complex structure.} Let us assume
$M^{2n}$ to be homogeneous space $G/H$ of positive Euler
characteristic with canonical action of a maximal torus and endowed
with an invariant almost complex structure $J$. All fixed points
have sign $+1$ and taking into account Theorem~\ref{weights},
Proposition~\ref{coeffzero1} gives that the weights at the fixed
points have to satisfy the following relations.

\begin{cor}
For any $\omega$ with $0\leqslant \|\omega\| \leqslant (n-1)$ where
$2n = \dim G/H$ we have that
\begin{equation}
\sum _{\rw \in W_{G}/W_{H}} \rw \Big( \frac{f_\omega(t_1, \ldots,
t_n)}{t_1 \cdots t_n}\Big ) = 0 \; ,
\end{equation}
where $t_j = \langle \Lambda _{j}, {\bf x}\rangle$ and $\Lambda
_{j}, \; 1\leqslant j\leqslant n$, are the weights  at the fixed
point $e \cdot H$.
\end{cor}

In the same way,  Theorem~\ref{sthm} implies that

\begin{thm}\label{s}
For $M^{2n}=G/H$ and $t_j=\langle  \Lambda_j,{\bf x} \rangle$, where
$\langle  \Lambda_j,{\bf x} \rangle=\sum\limits_{l=1}^k \Lambda_j^l
x_l, \; {\bf x}=(x_1,\ldots,x_k), \; k=\rk G=\rk H$, we have
\begin{equation}
s_{\omega}(\tau (M^{2n}))=\sum _{\rw \in W_{G}/W_{H}}\rw \Big(
\frac{f_\omega(t_1, \ldots, t_n)}{t_1 \cdots t_n}\Big)
\end{equation}
for any $\omega$ such that $\|\omega\|=n$.
\end{thm}

\begin{ex}\label{euler}
\[
s_{(n,0,\ldots,0)}(G/H, J)= \| W_{G}/W_{H}\| = \chi (G/H)
\]
and, therefore, $s_{(n,0,\ldots ,0)}(G/H, J)$ does not depend on
invariant almost complex structure $J$.
\end{ex}

\begin{cor}\label{sn}
\[
s_{(0,\ldots,0,1)}(G/H, J)= s_n(G/H, J))= \sum _{\rw \in
W_{G}/W_{H}} \rw \Big( \frac{\sum\limits_{j=1}^n t_j^n}{t_1 \cdots
t_n}\Big) \; .
\]
\end{cor}

\begin{ex}In the case $\mathbb{C}P^n=G/H$ where $G=U(n+1),\;H=U(1)\times
U(n)$ we have action of $T^{n+1}$ and related to the standard
complex structure the weights are given with $\langle\Lambda_j,{\bf
x}\rangle = x_j-x_{n+1}, \;j=1,\ldots, n$ and
$W_{G}/W_{H}=\mathbb{Z}_{n+1}$ is cyclic group. So
\begin{equation}
s_n(\mathbb{C}P^n) = \sum_{i=1}^{n+1}\frac{\sum\limits_{j\neq i}
(x_i-x_j)^n}{\prod\limits_{j\neq i}(x_i-x_j)}= n+1 \; .
\end{equation}
\end{ex}

\begin{ex}
Let us consider Grassmann manifold $G_{q+2,2}=G/H$ where
$G=U(q+2),\;H=U(q)\times U(2)$. We have here the canonical action of
the torus $T^{q+2}$. The weights for this action at identity point
related to the standard complex structure are given with $\langle
\Lambda _{ij}, {\bf x}\rangle = x_i - x_j$, where $1\leqslant i\leqslant q$, 
$j=q+1,q+2$. There are $\|W_{U(q+2)}/W_{U(2)\times
U(q)}\|=\frac{(q+2)(q+1)}{2}$ fixed points for this action.
Therefore
\begin{equation}\label{gn2}
s_{2q}(G_{q+2,2}) = \sum _{\rw \in W_{U(q+2)}/W_{U(2)\times U(q)}}
\rw \Big( \frac{\sum\limits_{i=1}^{q}\big(
(x_i-x_{q+1})^{2q}+(x_i-x_{q+2})^{2q}\big)
}{\prod\limits_{i=1}^{q}(x_i-x_{q+1})(x_i-x_{q+2})}\Big) \; .
\end{equation}
The action of the group  $W_{U(q+2)}/W_{U(2)\times U(q)}$ on the
weights at the identity  point in formula~\eqref{gn2} is given by
the permutations between the coordinates $x_1,\ldots ,x_q$ and coordinates $x_{q+1},x_{q+2}$. The explicit description of non trivial such permutations is as follows: 
\begin{align*}
\rw _{k,q+1}(k) &= q+1, \; \rw _{k,q+1}(q+1) =k, \; \mbox{where} \;1\leqslant
k\leqslant q \; ,\\
\rw _{k,q+2}(k) &= q+2, \; \; \rw _{k,q+2}(q+2)\; = k, \; \mbox{where} \;
1\leqslant k\leqslant q \; ,
\end{align*}
\[ \rw _{k,l}(q+1)=k, \;  \rw _{k,l}(k) = q+1, \; \rw _{k,l}(q+2)=l, \; \rw
_{k,l}(l) =q+2 \;\; \mbox{for}\; 1\leqslant k\leqslant q-1, \;
k+1\leqslant l\leqslant q \; . 
\]
As we remarked before (see Remark~\ref{int}), the expression on the
right hand side in~\eqref{gn2} is an integer number, so we can get a value
for $s_{2q}$ by choosing the appropriate values for the vector
$(x_1,\ldots ,x_{q+2})$. For example, if we take $q=2$ and
$(x_1,x_2,x_3,x_4) = (1,2,3,4)$ the straightforward application of
formula~\eqref{gn2} will give that $s_{4}(G_{4,2}) = -20$.
\end{ex}

\begin{ex}In the case $G_{q+l,l}=G/H$, where $G=U(q+l),\; H=U(q)\times
U(l)$ we have
\begin{equation}
s_{lq}(G_{q+l,l}) = \sum_{\sigma \in S_{q+l}/(S_q \times S_l)}\sigma
\Big( \frac{\sum (x_i-x_j)^{lq}}{\prod(x_i-x_j)}\Big) \; ,
\end{equation}
where $1 \leqslant i \leqslant q, \; (q+1) \leqslant j \leqslant
(q+l)$ and $S_{q+l}$ is the symmetric group.
\end{ex}

We consider later, in the Section~\ref{app}, the case of this
Grassmann manifold in more details.

\subsubsection{Chern numbers.}
We want to deduce an explicit relations between cohomology
characteristic numbers $s_{\omega}$ and classical Chern numbers for
an invariant almost complex structure on $G/H$.

\begin{prop}\label{orbit}
The number $s_{\omega}(\tau (M^{2n}))$, where $\omega=(i_1, \ldots,
i_n), \; \|\omega\|=n$, is the characteristic number that
corresponds to the characteristic class given by the orbit of the
monomial
\[
(u_{1}\cdots u_{i_1})(u_{i_1+1}^2\cdots u_{i_1+i_2}^2)\cdots
(u_{i_1+\ldots +i_{n-1}+1}^{n}\cdots u_{i_1+\ldots +i_n}^{n}) \; .
\]
\end{prop}

\begin{rem}
Let $\xi=(j_1,\ldots, j_n)$ and ${\bf u}^\xi=u_{1}^{j_1}\cdots
u_{n}^{j_n}$. The orbit of the monomial ${\bf u}^\xi$ is defined
with
\[
O({\bf u}^\xi) = \sum {\bf u}^{\xi'} \; ,
\]
where the sum is over the orbit $\{\xi'=\sigma \xi, \; \sigma \in
S_n \}$ of the vector $\xi \in \mathbb{Z}^n$ under the symmetric
group $S_{n}$ acting by permutations of coordinates of $\xi$.
\end{rem}

\begin{ex}\label{ex1}
If we take $\omega = (n,0,\ldots ,0)$ we need to compute the
coefficient for $a_{1}^{n}$ and it is given as an orbit
$O(u_{1}\cdots u_{n})$ what is the  elementary symmetric function
$\sigma _{n}$. If we take $\omega =(0,\ldots,0,1)$ then we should
compute the coefficient for $a_{n}$ and it is given with
$O(u_{1}^{n})=\sum _{j=1}^{n} u_{j}^{n}$, what is Newton polynomial.
\end{ex}

It is well known fact from the algebra of symmetric functions that
the orbits of monomials give the additive basis for the algebra of
symmetric functions. Therefore, any orbit of monomial can be
expressed through elementary symmetric functions and vice versa. It
gives the expressions  for the characteristic numbers $s_{\omega}$ in
terms of Chern characteristic numbers $c^\omega=c_1^{i_1}\cdots
c_n^{i_n}$ for an almost complex homogeneous space $(G/H, J)$.

\begin{thm}\label{sc}
Let $\omega =(i_1,\ldots, i_n), \; \|\omega\|=n$, and assume that
the orbit of the monomial
\[
(u_{1}\cdots u_{i_1})(u_{i_1+1}^2\cdots u_{i_1+i_2}^2)\cdots
(u_{i_1+\ldots +i_{n-1}+1}^{n}\cdots u_{i_1+\ldots +i_n}^{n})
\]
 is expressed through the elementary symmetric functions as
\begin{equation}\label{change}
O((u_{1}\cdots u_{i_1})(u_{i_1+1}^2\cdots u_{i_1+i_2}^2)\cdots
(u_{i_1+\ldots +i_{n-1}+1}^{n}\cdots u_{i_1+\ldots +i_n}^{n}))=
\end{equation}
\[
= \sum _{\| \xi \| = n} \beta_{\omega \xi}\sigma _{1}^{l_1}\cdots \sigma
_{n}^{l_n}
\]
for some $\beta_{\omega \xi}\in \Z$ and $\| \xi\| = \sum _{j=1}^{n}j\cdot
l_j$, where $\xi=(l_1, \ldots, l_n)$. Then  it holds
\begin{equation}\label{nocohomchern}
s_{\omega}(G/H, J) = \sum _{\rw \in W_{G}/W_{H}}\rw \Big(
\frac{f_\omega(t_1, \ldots, t_n)}{t_1 \dots t_n}\Big)= \sum _{\| \xi
\| = n} \beta_{\omega \xi}c_{1}^{l_1}\cdots c_{n}^{l_n} \; ,
\end{equation}
where $c_{i}$ are the Chern classes for the tangent bundle of
$(G/H,J)$.
\end{thm}

\begin{rem}
Let $p(n)$ denote the number of  partitions of the number $n$. By
varying $\omega$, the equation~\eqref{nocohomchern} gives the system
of $p(n)$ linear equations in Chern numbers whose determinant is,
by~\eqref{change},  non-zero. Therefore, it provides the explicit
formulas for the computation of Chern numbers.
\end{rem}

\begin{rem}
We want to point that relation~\eqref{nocohomchern} in 
Theorem~\ref{sc} together with Theorem~\ref{s} proves that the Chern
numbers for $(G/H, J)$ can be computed without having  any
information on cohomology for $G/H$.
\end{rem}

\begin{ex}\label{ex4}
We provide the  direct application of Theorem~\ref{sc} following
Example~\ref{ex1}. It is straightforward to see  that
$s_{(n,0,\ldots ,0)}(G/H) = c_{n}(G/H)$ for any invariant almost
complex structure. This together with Example~\ref{euler} gives that
$c_{n}(G/H) = \chi (G/H)$.
\end{ex}

We want to add that it is given in~\cite{ms} a description of the
numbers $s_{I}$ that correspond to our characteristic numbers
$s_{\omega}$, but the numerations $I$ and $\omega$ are  different.
To the partition $i\in I$ correspond the n-tuple $\omega
=(i_1,\ldots ,i_n)$ such that $i_k$ is equal to the number of
appearances of the number $k$ in the partition $i$.

\section{On equivariant stable complex homogeneous spaces of
positive Euler characteristic.}\label{eqsthom}

Assume we are given on $G/H$ a stable complex structure $c_{\tau}$
which is equivariant under the canonical action $\theta$ of the
maximal torus $T$. If $p=gH$ is the fixed point for an action
$\theta$, by  Section 2 and Section 3, we see  that $T_{gH}(G/H)=\gg
^{\C}_{\rw (\alpha _{1})}\oplus \ldots \oplus\gg ^{\C}_{\rw (\alpha
_{n})}$, where $\rw \in W_{G}/W_{H}$ and $\alpha _{1},\ldots ,\alpha
_{n}$ are the complementary roots for $G$ related to $H$. The
following statement is the direct consequence of 
Theorem~\ref{weights2} and Lemma~\ref{wfixed}.

\begin{cor}\label{stableweights}
Let $\alpha _{1},\ldots ,\alpha _{n}$ be the set of complementary
roots for $G$ related to $H$ . The set of weights of an action
$\theta$ at the fixed points related to an arbitrary equivariant
stable complex structure $c_{\tau}$ is of the form
\begin{equation}\label{stablevectors}
\{ a_{1}(\rw )\cdot \rw (\alpha _{1}),\ldots ,a_{n}(\rw )\cdot
\rw (\alpha _{n}) \} \; ,
\end{equation}
where  $\rw \in W_{G}/W_{H}$ and $a_{i}(\rw )=\pm 1$ for $1\leqslant
i\leqslant n$.
\end{cor}

For the signs at the fixed points  (which we identify with $\rw \in
W_{G}/W_{H}$) using  Theorem~\ref{weights2} we obtain the following.

\begin{cor}\label{stablesigns}
Assume that $G/H$ admits an invariant almost complex structure $J$
defined by the complementary roots $\alpha _{1},\ldots,\alpha_{n}$.
Let $c_{\tau}$ be the $\theta$-equivariant stable complex structure
with the set of weights $\{ a_{1}(\rw )\cdot \rw (\alpha
_{1}),\ldots ,a_{n}(\rw )\cdot \rw (\alpha _{n}) \}$, $\rw \in
W_{G}/W_{H}$ at the fixed points. The signs at the fixed points are
given with
\begin{equation}\label{csign}
\sg (\rw ) = \epsilon \cdot \prod _{i=1}^{n}a_{i}(\rw ), \;\; \rw\in W_{G}/W_{H} \; ,
\end{equation}
where $\epsilon =\pm 1$ depending if $J$ and $c_{\tau}$ define the
same orientation on $M^{2n}$ or not.
\end{cor}

This implies the following consequence of  Proposition~\ref{neccond}.

\begin{cor}\label{stablecoeffzero}
Assume that $G/H$ admits an invariant almost complex structure $J$
defined by the complementary roots $\alpha _{1},\ldots,\alpha_{n}$.
Let $c_{\tau}$ be the $\theta$-equivariant stable complex structure
with the set of weights $\{ a_{1}(\rw )\cdot \rw (\alpha
_{1}),\ldots ,a_{n}(\rw )\cdot \rw (\alpha _{n}) \}$, $\rw \in
W_{G}/W_{H}$ at the fixed points. Then
\begin{itemize}
\item for any $\omega$ with  $0\leqslant \|\omega\| \leqslant (n-1)$ we have that
\begin{equation}
\sum _{\rw \in W_{G}/W_{H}}\frac{f_{\omega}(a_{1}(\rw )\rw (\alpha _{1}),\ldots ,
a_{n}(\rw )\rw (\alpha _{n}))}{\rw (\alpha_1) \cdots \rw (\alpha_n)} = 0 \; ;
\end{equation}
\item for any $\|\omega\|=n$
\begin{equation}
\sum _{\rw \in W_{G}/W_{H}}
\frac{f_\omega(a_{1}(\rw )\rw (\alpha _{1}),\ldots, a_{n}(\rw )\rw (\alpha _{n}))}
{\rw (\alpha _{1}) \cdots \rw (\alpha_{n})}
\end{equation}
is an integer number.
\end{itemize}
\end{cor}

\begin{rem}
Corollary~\ref{stablecoeffzero} gives the strong constraints on
the numbers $a_{i}(\rw )$ that appear as the "coefficients" (related
to $J$) of the weights of the $\theta$-equivariant stable complex
structure. In that way, it provides information which  vectors $\{
a_{1}(\rw )\cdot \rw (\alpha _{1}),\ldots ,a_{n}(\rw )\cdot \rw
(\alpha _{n}) \}$, $\rw \in W_{G}/W_{H}$ can not be realized as the
weight vectors of some $\theta$ equivariant stable complex
structure.
\end{rem}

\begin{ex}\label{cps}
Consider complex projective space $\C P^{3}$ and let $c_{\tau}$ be
a stable complex structure on $\C P^{3}$ equivariant under the
canonical action of the torus $T^4$. By
Corollary~\ref{stableweights}, the weights for $c_{\tau}$ at the
fixed points are $a_{1}(\rw )\cdot \rw (x_1-x_4), a_{2}(\rw )\cdot
\rw (x_2-x_4), a_{3}(\rw )\cdot \rw (x_3-x_4)$, where $\rw \in
W_{U(4)}/W_{U(3)} = \Z _{4}$. Corollary~\ref{stablecoeffzero}
implies that the coefficients in $t$ and $t^2$ in the polynomial
\begin{equation}\label{pol}
\prod_{1\leq i<j\leq 4}(x_i-x_j)\cdot \sum_{\rw \in \Z_{4}}
\frac{f(ta_{1}(\rw )\rw (x_1-x_4))f(ta_{2}(\rw )\rw
(x_2-x_4))f(ta_{3}(\rw )\rw (x_3-x_4))}{ \rw (x_1-x_4)\rw
(x_2-x_4)\rw (x_3-x_4)} \; ,
\end{equation}
where $f(t)=1+a_1t+a_{2}t^2+a_{3}t^3$, have to be zero. The
coefficient in $t$ for~\eqref{pol} is a polynomial
$P(x_1,x_2,x_3,x_4)$ of degree $4$ whose coefficients are some
linear combinations of the numbers $a_{i}(\rw )$. Some of them are
as follows
\[
x_1x_4^3:a_{2}(2) - a_{3}(3), \;\; x_1^3x_4: a_{1}(3)-a_{1}(2), \;\;
x_1^3x_2: a_{1}(0)-a_{1}(3), \;\; x_2^3x_4: a_{2}(1)-a_{2}(3),
\]
\[
x_3x_4^3: a_{1}(1)-a_{2}(2),\;\; x_3^3x_4: a_{3}(2)-a_{3}(1),\;\;
x_1x_3^3: a_{3}(0)-a_{3}(2) \; ,
\]
what implies that
\begin{equation}\label{arelations}
a_{1}(0)=a_{1}(3)=a_{1}(2),\; a_{2}(0) = a_{2}(1) = a_{2}(3) \; ,
\end{equation}
\[ a_{3}(0) = a_{3}(1)=a_{3}(2), a_{1}(1)= a_{2}(2)=a_{3}(3) \; .
\]
The direct computation shows that the requirement
$P(x_1,x_2,x_3,x_4)\equiv 0$ is equivalent to the
relations~\eqref{arelations} and also that the same relations  give
that the coefficient in $t^2$ in the  polynomial~\eqref{pol}  will
be zero.

Therefore, the weights for $c_{\tau}$  at the fixed points are
completely determined by the values for $a_{1}(0), a_{2}(0),
a_{3}(0), a_{1}(1)$. By Corollary~\ref{stablesigns} the signs at the
fixed points are also determined by these  values up to factor
$\epsilon = \pm 1$ depending if $c_{\tau}$ and standard canonical structure
on $\C P^{3}$ give the same orientation or not.

In turns out that in this case  for all possible values $\pm 1$ for
$a_{1}(0), a_{2}(0), a_{3}(0), a_{1}(1)$ the corresponding
vectors~\eqref{stablevectors} can be realized as the weights vectors
of the stable complex structures. In order to verify this we look at
the decomposition  $T(\C P^{3})\oplus \C \cong \hat {\eta}\oplus
\hat{\eta}\oplus \hat{\eta}\oplus \hat{\eta}$, where $\hat{\eta}$ is
the conjugate to the Hopf bundle over $\C P^{3}$. Using this
decomposition  we can get the stable complex structures on $T(\C
P^{3})$ by choosing the complex structures on each $\hat{\eta}$.
The weight vector for any such stable complex structure is
determined by the corresponding choices for the values of $a_{1}(0),
a_{2}(0), a_{3}(0), a_{1}(1)$. For example if we take
$a_{1}(0)=a_{2}(0)=a_{3}(0)=a_{1}(1)=1$ we get the weights for the
standard complex structure which we can get from the stable complex
structure  $\hat {\eta}\oplus \hat{\eta}\oplus \hat{\eta}\oplus
\hat{\eta}$.  

If we take $a_{1}(0)=a_{2}(0)=a_{3}(0)=1$,
$a_{1}(1)=-1$ the corresponding weights come from the stable complex
structure $\hat{\eta}\oplus \hat{\eta}\oplus \hat{\eta}\oplus \eta$.
It follows  that this stable complex structure and the standard complex structure 
define an opposite orientations on $\C P^3$ and, therefore, $\epsilon=-1$.   
Usign~\eqref{csign} we obtain that related to this stable complex structure the signs at the fixed
points are $\sg(0)=-1$, $\sg(1)=\sg(2)=\sg(3)=1$. The weights at the fixed points are:
\[
(0) : x_1-x_4,\; x_2-x_4,\; x_3-x_4;\;\; (1) : x_1-x_4,\; x_2-x_1,\; x_3-x_1\;;
\]
\[
(2) : x_1-x_2,\; x_2-x_4,\; x_3-x_2;\;\; (3) : x_1-x_3,\; x_2-x_3,\; x_3-x_4\;.
\]
By Example~\ref{ex3} we obtain that the number $s_{3}$ for $\C P^3$ related to this stable complex 
structure is equal to $-2$. In that way it shows that $\C P^3$ with  this non-standard stable complex structure realizes multiplicative generator in complex cobordism ring of dimension $6$. 

The relations~\eqref{arelations}  also give an examples of  the
vectors~\eqref{stablevectors} that can not be realized as the weight
vectors of the stable complex structures on $\C  P^{3}$ equivariant
under the canonical torus action.
\end{ex}

\begin{ex}
Using Proposition~\ref{stableweights}, we can also find an examples
of the vectors~\eqref{stablevectors} on the flag manifold $U(3)/T^3$
or Grassmann manifold $G_{4,2}$ that can {\it not} be realized as
the weight vectors of the stable complex structures.

In a case of $U(3)/T^3$ any stable complex structure  has the
weights of the form $a_1(\sigma)\sigma (x_1-x_2),a_2(\sigma)\sigma
(x_1-x_3),a_3(\sigma)\sigma (x_2-x_3)$, where $\sigma \in S_{3}$.
The relation~\eqref{secondstable} gives that for $k=1$  the numbers
$a_{i}(\sigma )$ have to satisfy the following relation
\[
\sum_{\sigma \in S_{3}}\frac{a_1(\sigma)\sigma
(x_1-x_2)+a_2(\sigma)\sigma (x_1-x_3)+a_3(\sigma)\sigma
(x_2-x_3)}{\sigma (x_1-x_2)\sigma (x_1-x_3)\sigma (x_2-x_3)} = 0 \;
.
\]
It is equivalent to
\[
a_{1}(123)+a_{2}(123)+a_1(213)-a_3(213)+a_2(321)+a_3(321) \]
\[-a_1(132)-a_2(132)-a_2(231)-a_3(231)-a_1(312)+a_3(312)= 0 \; ,
\]
\[
-a_{1}(123)+a_{3}(123)-a_1(213)-a_2(213)+a_1(321)-a_3(321) \]
\[+a_2(132)-a_3(132)+a_1(231)+a_2(231)-a_2(312)-a_3(312)= 0 \; .
\]
It follows from these relations that the
vector~\eqref{stablevectors} determined with $a_{1}(123)=-1$ and
$a_{i}(\sigma)=1$ for all the others $1\leqslant i\leqslant 3$ and
$\sigma \in S_{3}$, can not be obtained as the weight vector of some
stable complex structure on $U(3)/T^3$ equivariant under the
canonical action of the maximal torus.

Using the same argument we can also conclude that for the
Grassmannian $G_{4,2}$ the vector determined with $a_{1}(1234)=-1$
and $a_{i}(\sigma )=1$ for all the others $1\leqslant i\leqslant 4$
and $\sigma \in S_{4}/S_{2}\times S_{2}$ can not be realized as  the
weight vector of the stable complex structure equivariant under the
canonical action of $T^4$.
\end{ex}

\section{Some applications.}\label{app}

\subsection{Flag manifolds $U(n)/T^n$.}
We consider invariant complex structure on $U(n)/T^n$.
Recall~\cite{Adams} that the
Weyl group $W_{U(n)}$ is the symmetric group and it permutes the
coordinates $x_1,\ldots ,x_n$ on Lie algebra $\TT ^{n}$ for $T^{n}$.
The canonical action of the torus $T^{n}$ on this manifold has $\|
W_{U(n)}\| = \chi (U(n)/T^n) = n!$ fixed points and its weights at
identity point are given by the roots of U(n).

We first consider the case $n=3$ and apply our results to explicitly
compute cobordism class and Chern numbers for $U(3)/T^3$. The roots
for U(3) are $x_1-x_2$, $x_1-x_3$ and $x_2-x_3$. Therefore the
cobordism class for $U(3)/T^3$ is given as the coefficient for $t^3$
in the polynomial
\[
[U(3)/T^3] = \sum _{\sigma \in S_3}\sigma
\Big(\frac{f(t(x_1-x_2))f(t(x_1-x_3))f(t(x_2-x_3))}{(x_1-x_2)
(x_1-x_3)(x_2-x_3)}\Big) \; ,
\]
where $f(t)=1+a_1t+a_2t^2+a_3t^3$, what implies
\[
[U(3)/T^3] = 6(a_1^3 +a_1a_2-a_3) \; .
\]
This gives that the characteristic numbers $s_\omega$ for $U(3)/T^3$
are
\[
s_{(3,0,0)}=6, \quad s_{(1,1,0)} = 6, \quad s_{(0,0,1)}=-6 \; .
\]
By  Theorem~\ref{sc} we have the following relations between
characteristic numbers $s_\omega$ and  Chern numbers $c^\omega$
\[
c_3 = 6, \; c_1c_2 - 3c_3=6, \; c_1^3-3c_1c_2+3c_3=-6, \; \;
\mbox{what gives}\; \; c_1c_2=24, \; c_1^3=48 \; .
\]
To simplify the notations we take further $\Delta_n =
\prod\limits_{1\leqslant i<j\leqslant n}(x_i-x_j)$.

\begin{thm}\label{thm}
The Chern-Dold character of the toric genus for the flag manifold
$U(n)/T^n$ is given by the formula:
\begin{equation}\label{40}
ch_{U}\Phi(U(n)/T^n) = \frac{1}{\Delta_n}\sum _{\sigma \in S_{n}}
\sg (\sigma) \sigma \Big(\prod\limits_{1\leqslant i<j\leqslant
n}f(x_i-x_j)\Big) \; ,
\end{equation}
where $f(t) = 1+\sum\limits_{i\geqslant 1}a_{i}t^{i}$ and $\sg
(\sigma )$ is the sign of the permutation $\sigma$.
\end{thm}

\subsubsection{Using of divided difference operators.}
Consider the ring of the symmetric polynomials $\Sym_n \subset
\mathbb{Z}[x_1, \ldots, x_n]$. There is a linear operator (see
\cite{Macdonald-95})
\[ L : \mathbb{Z}[x_1, \ldots, x_n] \longrightarrow \Sym_n  : \;
 L{\bf x}^\xi = \frac{1}{\Delta_n }
\sum _{\sigma \in S_{n}} \sg (\sigma) \sigma ({\bf x}^\xi) \; , \] where
$\xi=(j_1, \ldots, j_n)$ and ${\bf x}^\xi = x_1^{j_1} \cdots
x_n^{j_n}$.

It follows from the definition of Schur polynomials $\Sh(x_1,
\ldots, x_n)$ where $\lambda=(\lambda_1 \geqslant \lambda_2
\geqslant \cdots \geqslant \lambda_n \geqslant 0)$ (see
\cite{Macdonald-95}), that
\[ L{\bf x}^{\lambda+\delta} = \Sh(x_1,\ldots, x_n) \; , \]
where $\delta=(n-1,n-2,\ldots,1,0)$ and $L{\bf x}^\delta=1$.
Moreover, the operator $L$ have the following properties:
\begin{itemize}
\item $L{\bf x}^\xi=0$, if $j_1 \geqslant j_2 \geqslant \cdots
\geqslant j_n\geqslant 0$ and $\xi \neq \lambda+\delta$ for some
$\lambda=(\lambda_1 \geqslant \lambda_2 \geqslant \cdots \geqslant
\lambda_n \geqslant 0)$;
\item $L{\bf x}^\xi=\sg (\sigma)L\sigma({\bf x}^\xi)$, where $\xi=(j_1,
\ldots, j_n)$ and $\sigma\xi=\xi'$, where $\sigma \in S_n$ and
$\xi'=(j'_1 \geqslant j'_2 \geqslant \cdots \geqslant j'_n\geqslant
0)$;
\item $L$ is a homomorphism of $\Sym_n$-modules.
\end{itemize}
We have
\begin{equation}\label{P}
\prod\limits_{1\leqslant i<j\leqslant n}f(t(x_i-x_j))=
1+\sum_{|\xi|>0} P_\xi(a_1, \ldots, a_n, \ldots)t^{|\xi|}{\bf x}^\xi
\; ,
\end{equation}
 where
$|\xi|=\sum\limits_{q=1}^n j_q$.

\begin{cor} \label{L}
The Chern-Dold character of the toric genus for the flag manifold
$U(n)/T^n$ is given by the formula:
\[ ch_{U}\Phi(U(n)/T^n) = \sum_{|\lambda|\geqslant 0}\Big(
\sum_{\sigma\in S_n}\sg (\sigma) P_{\sigma(\lambda+\delta)}(a_1,
\ldots, a_n, \ldots)\Big)\Sh(x_1, \ldots, x_n) \,, \] where
$\delta=(n-1,n-2,\ldots,1,0),\; \lambda=(\lambda_1 \geqslant
\lambda_2 \geqslant \cdots \geqslant \lambda_n \geqslant 0)$. In
particular,
\begin{equation}\label{CL}
[U(n)/T^n] = \sum_{\sigma\in S_n}\sg (\sigma) P_{\sigma\delta}(a_1,
\ldots, a_n, \ldots)\,.
\end{equation}
\end{cor}
\begin{proof}\label{delta}
Set $m=\frac{n(n-1)}{2}$. From  Theorem \ref{thm} and the formula
\eqref{P} we obtain:
\[ ch_{U}\Phi(U(n)/T^n) =  \sum_{|\xi|\geqslant m} P_\xi L{\bf x}^\xi\,. \]
The first property of the operator $L$ gives that for any $\xi$ we
will have $L{\bf x}^{\xi}=0$, whenever ${\bf x}^{\xi}\neq \sigma
({\bf x}^{\lambda+\delta})$ for some $\sigma \in S_n$ and
$\lambda=(\lambda_1 \geqslant \lambda_2 \geqslant \cdots \geqslant
\lambda_n \geqslant 0)$. The second property gives that
$L\sigma({\bf x}^{\lambda+\delta})= \sg (\sigma)\Sh(x_1, \ldots,
x_n)$.
\end{proof}

\begin{rem}
1. In the case $n=2$ this corollary gives the result of Example
\ref{n=2}.\\
2. As we will show in Corollary~\ref{cor8} below, polynomials
$P_{\sigma\delta}$ in the formula~\eqref{CL} appears to be
polynomials only in variables $a_1,\ldots, a_{2n-3}$.
\end{rem}

Set by definition
\[ 1+\sum_{|\xi|>0} \sigma^{-1}(P_\xi)t^{|\xi|}{\bf x}^\xi =
\sigma\Big( \prod\limits_{1\leqslant i<j\leqslant n}f(t(x_i-x_j))
\Big), \] where $\sigma \in S_n$ on the right acts by the permutation
of variables $x_1, \ldots, x_n$. Directly from the definition we
have
\[ 1+\sum_{|\xi|>0} \sigma^{-1} (P_\xi)t^{|\xi|}{\bf x}^\xi = 1+\sum_{|\xi|>0}P_\xi t^{|\xi|}\sigma({\bf
x}^\xi). \] Therefore
\begin{equation}\label{sigma}
\sigma (P_\xi) = P_{\sigma \xi}\;.
\end{equation}
Together with Corollary~\ref{L} the formula~\eqref{sigma} implies the following Theorem.

\begin{thm}\label{t-chi}
\begin{equation}\label{chi}
[U(n)/T^n] = \Big(\sum_{\sigma\in S_n}\sg (\sigma)\sigma\Big)
P_{\delta}(a_1, \ldots, a_n, \ldots)\;.
\end{equation}
\end{thm}

\begin{cor}
\[ \sigma[U(n)/T^n]=\sg (\sigma)[ U(n)/T^n]\;. \]
\end{cor}
\begin{proof}
\[
\sigma[U(n)/T^n]=\Big(\sum_{\tilde{\sigma}\in S_n}\sg (\tilde{\sigma})\sigma\tilde{\sigma}\Big)P_{\delta}(a_1,\ldots ,a_n,\ldots)=
\]
\[
=\Big(\sum_{\bar{\sigma}\in S_n}\sg (\sigma^{-1}\bar{\sigma})\bar{\sigma}\Big)P_{\delta}(a_1,\ldots ,a_n,\ldots)=
\sg (\sigma^{-1})\Big(\sum_{\bar{\sigma}\in S_n}\sg (\bar{\sigma})\bar{\sigma}\Big)P_{\delta}(a_1,\ldots ,a_n,\ldots) \; .
\]
Since $\sg (\sigma) = \sg (\sigma^{-1})$, the formula follows.
\end{proof}

\begin{ex}
The direct computation gives that for $n=3$ the polynomials $P_{\sigma\delta}(a_1,a_2,a_3)=\sigma(P_{\delta})$ from~\eqref{P}, where $\delta = (2,1,0)$ and $\sigma \in S_3$ are
\[
P_{\delta}=a_1^3-a_1a_2-3a_3\;,
\]
\[
(12)P_{\delta}=(13)P_{\delta}=-P_{\delta}\;, 
\]
\[
(23)P_{\delta}=-(a_1^3+5a_1a_2+3a_3)\;.
\]
The action of the transpositions implies that the rest two permutations act as
\[
(312)P_{\delta}=P_{\delta},\quad (231)P_{\delta}=-(23)P_{\delta} \;.
\]
Using Corollary~\ref{L} we obtain the cobordism class $[U(3)/T^3]=
6(a_1^3+a_1a_2-a_3)$. The above also gives that the symmetric group $S_{3}$ acts
non-trivially on $P_{\delta}(a_1,a_2,a_3)$ and
$\sum\limits_{\sigma\in S_3}\sigma P_\delta(a_1,a_2,a_3)=0$.
\end{ex}

We have
\[  \prod\limits_{1\leqslant i<j\leqslant n}f(t(x_i-x_j)) \equiv
\prod\limits_{1\leqslant i<j\leqslant n}f(t(x_i+x_j))\mod 2\;. \]
Using that $ \prod\limits_{1\leqslant i<j\leqslant n}f(t(x_i+x_j))$
is the symmetric series of variables $x_1,\ldots,x_n$ we obtain
\[ \sigma(P_\xi) \equiv P_\xi \mod 2\;, \]
for any $\sigma\in S_n$. Thus Theorem \ref{t-chi} implies the
following:

\begin{cor}
All Chern numbers of the manifold $U(n)/T^n,\; n \geqslant 2$, are
even.
\end{cor}

\begin{rem}
Cobordism class $[U(3)/T^3]$ gives nonzero element in
$\Omega_U^{-6}\otimes \mathbb{Z}/2$ because\\ $s_{3}(U(3)/T^3)\equiv
2 \mod 4$.
\end{rem}

The characteristic number $s_{m}$ for $U(n)/T^{n}$ is given as
\begin{equation} \label{sm}
s_{m}(U(n)/T^n) = \sum _{1\leqslant i<j\leqslant n} L(x_i-x_j)^{m} \;
.
\end{equation}
Corollary \ref{L} implies the following:

\begin{cor}
$s_{1}(U(2)/T^2)=2; \; s_{3}(U(3)/T^3)= -6$\; and
\begin{equation}
s_{m}(U(n)/T^n) = 0 \; ,
\end{equation}
where $m=\frac{n(n-1)}{2}$ and $n>3$.
\end{cor}

We can push up this further. Denote by $(u_1,\ldots,u_m) = \big(
(x_i-x_j),\; i<j \big)$, where $m=\frac{n(n-1)}{2}$. Then for
$\omega = (i_1,\ldots ,i_m),\, \|\omega\|=m$ we have that
\[ O\Big( (u_1 \cdots u_{i_1})(u_{i_1+1}^2
\cdots u_{i_1+i_2}^2)\cdots (u_{i_{1}+\ldots +i_{m-1}+1}^m\cdots u_{i_{1}+\ldots +i_{m}}^m)\Big) =
\sum_{|\xi|=m}\alpha_{\omega,\xi}{\bf x}^\xi\;. \] This implies that
\begin{equation}
s_{\omega}(U(n)/T^n) = \sum_{|\xi|=m}\alpha_{\omega,\xi}L{\bf x}^\xi = \sum_{\sigma\in S_n}\sg (\sigma)
\alpha_{\omega,\sigma\delta}\,.
\end{equation}
Therefore, if $\xi=(j_1, \ldots, j_n)$, then
$\max\limits_{p_1,\ldots,p_s}(j_{p_1}+ \cdots+ j_{p_s}) = s\Big(
n-\frac{s+1}{2} \Big), \; 1\leqslant s \leqslant n$. In particular,
it holds that $\max\limits_{p_1,p_2}(j_{p_1}+j_{p_2}) = 2n-3$.

\begin{cor} \label{cor8}
Let $\omega = (i_1,\ldots ,i_m)$ such that $i_k\neq 0$ for some
$k>2n-3$, then
\[ s_{\omega}(U(n)/T^n) = 0 \; .\]
\end{cor}

If $\omega = (i_1,\ldots ,i_k), \,\|\omega\|=m$, does not satisfy
Corollary~\ref{cor8}, but $i_{k_1},\ldots ,i_{k_l} \neq 0$ for some
$k_1,\ldots ,k_l$ then we have that $k_{p}=2(n-1)-q_{p}$, for
$q_p\geqslant 1$, $1\leqslant p\leqslant l$. In this case we can say the following.

\begin{cor}\label{cor9}
If $n\geqslant 2l$ and $\sum\limits_{p=1}^{l}q_{p} < l(2l-1)$ then
\[ s_{\omega}(U(n)/T^n) = 0\; .\]
\end{cor}

\begin{rem}\label{trans}
From the second property of the operator $L$ we obtain that $L
\mathcal{P}(x_1, \ldots, x_n)=0$ for any series $\mathcal{P}(x_1,
\ldots, x_n)$, whenever $\sigma (\mathcal{P}(x_1, \ldots, x_n)) =
\varepsilon \mathcal{P}(x_1, \ldots, x_n)$ for a permutation $\sigma
\in S_n$, where $\varepsilon= \pm 1$ and $\varepsilon \cdot
\sg(\sigma)=-1$. This, in particular, gives that
$L(\mathcal{P}(x_1,\ldots ,x_n)+\sigma _{ij}(\mathcal{P}(x_1,\ldots
,x_n)))=0$ for any transposition $\sigma _{ij}$ of $x_i$ and $x_j$,
where $1\leqslant i<j\leqslant n$.
\end{rem}

Using Remark~\ref{trans} we can compute some more characteristic
numbers of the flag manifolds.

\begin{cor}
Let $n=4q$ or $4q+1$ and $\omega=(i_1, \ldots, i_m), \; \| \omega \|
= m$, where $i_{2l-1}=0$ for $l=1, \ldots, \frac{m}{2}$. Then
$s_{\omega}(U(n)/T^n) = 0$.
\end{cor}

Since $\sigma
_{12}((x_1-x_2)^{2l}\prod\limits_{\underset{(i,j)\neq(1,2)}{1\leqslant
i<j\leqslant n}}f(t(x_i-x_j)))=(x_1-x_2)^{2l}
\prod\limits_{\underset{(i,j)\neq(1,2)}{1\leqslant i<j\leqslant
n}}f(t(x_i-x_j))$ we have, also because of Remark~\ref{trans}, that
\[ L \Big(\prod\limits_{1\leqslant i<j\leqslant
n}f(t(x_i-x_j)\Big) = L \Big(\widetilde f(t(x_1-x_2))
\prod\limits_{\underset{(i,j)\neq(1,2)}{1\leqslant i<j\leqslant
n}}f(t(x_i-x_j))\Big)\; ,
\]
where $\widetilde f(t)=\sum\limits_{l\geqslant 1}a_{2l-1}t^{2l-1}$.
Using this property of $L$ once more we obtain

\begin{thm} \label{thm8}
For $n\geqslant 4$ the cobordism class for the flag manifold
$U(n)/T^n$ is given as the coefficient for $t^{\frac{n(n-1)}{2}}$ in
the series in $t$
\begin{equation} \label{tilde}
L \Big(\widetilde f(t(x_1-x_2))\widetilde f(t(x_{n-1}-x_n))
\prod\limits_{\underset{(i,j)\neq(1,2),(n-1,n)}{1\leqslant
i<j\leqslant n}}f(t(x_i-x_j))\Big)\; .
\end{equation}
\end{thm}

\begin{rem}
Corollary~\ref{cor9} implies that if $s_{\omega}\neq 0$ for some
$\omega =(i_1,\ldots ,i_m)$, then for some $1\leqslant l\leqslant
\frac{m}{2}$ it has to be $i_{2l-1}\neq 0$. The Theorem~\ref{thm8}
gives stronger results that, for $n\geqslant 4$ in polynomials
$P_{\sigma\delta}$ in~\eqref{CL} each monom contains the factor
$a_{2i_1-1}a_{2i_2-1}$.
\end{rem}

Theorem~\ref{thm8} provide a way for direct computation of the
number $s_{\omega}$, for $\omega =(i_1,\ldots ,i_m)$ such that $\|
\omega \| =2$, where $\| \omega \| = i_1+\ldots +i_m$.  For $n>5$ we
have that $s_{\omega}(U(n)/T^n) = 0$  for such $\omega$.  For $n=4$
and $n=5$ these numbers can be computed very straightforward as the
next example shows.

\begin{ex}
We provide
computation of the characteristic
number $s_{(1,0,0,0,1,0)}$ for $U(4)/T^4$. From the formula
(\ref{tilde}) we obtain immediately:
\begin{align*}
s_{(1,0,0,0,1,0)}(U(4)/T^4) &= L\Big( (x_1-x_2)(x_3-x_4)^5 +
(x_1-x_2)^5(x_3-x_4) \Big)=\\
&=10L \Big( (x_1-x_2)(x_3-x_4)(x_1^2x_2^2+x_3^2x_4^2)\Big)=\\
&=20L \Big( x_1^3x_2^2(x_3-x_4)+(x_1-x_2)x_3^3x_4^2)\Big)=\\
&=40L \Big( x_1^3x_2^2x_3+x_1x_3^3x_4^2\Big)=80\; .
\end{align*}
\end{ex}

\begin{rem}
We want to emphasize that the formula~\eqref{40} gives the
description of the cobordism classes of the flag manifolds in terms
of {\it divided difference operators}. The divided difference
operators are defined with (see \cite{BGG})
\[\partial_{ij}P(x_1,\ldots ,x_n)=\frac{1}{x_i-x_j}\Big(P(x_1,\ldots
,x_n)-\sigma _{ij}P(x_1,\ldots ,x_n)\Big),\] where $i<j$. Put
$\sigma _{i,i+1}=\sigma_i$, $\partial_{i,i+1}=\partial_i, \;
1\leqslant i \leqslant n-1$. We can write down operator $L$ as the
following composition (see \cite{F, Macdonald-91})
\[
L= (\partial_1 \partial_2 \cdots \partial_{n-1})(\partial_1 \partial_2
\cdots \partial_{n-2}) \cdots (\partial_1 \partial_2) \partial_1 \; .
\]

Denote by $\rw _0$ the permutation $(n,n-1, \ldots, 1)$. Wright down
a permutation $\rw \in S_n$ in the form $\rw = \rw _0 \sigma _{i_1}
\cdots \sigma _{i_p}$ and set $\nabla_{\rw} =
\partial_{i_p}\cdots
\partial_{i_1}$. It is natural to set $\nabla_{\rw _0} = I$  --- identity
operator. The space of operators $\nabla_{\rw}$ is dual to the space
of the Schubert polynomials
$\mathfrak{G}_{\rw}=\mathfrak{G}_{\rw}(x_1,\ldots ,x_n)$, since it
follows from their definition that $\mathfrak{G}_{\rw} =
\nabla_{\rw} {\bf x}^\delta$. Note that $\mathfrak{G}_{\rw _0}= {\bf
x}^\delta$. For the identity permutation $e=(1,2,\ldots, n)$ we have
$e=\rw _0 \cdot \rw _0^{-1}$. So $\nabla_e = L$ and
$\mathfrak{G}_{e} = \nabla_{e} {\bf x}^\delta = 1$.

Schubert polynomials were introduced in \cite{BGG} and in \cite{D}
in context of an arbitrary root systems. The main reference on
algebras of operators $\nabla_{\rw}$ and Schubert polynomials
$\mathfrak{G}_{\rw}$ is \cite{Macdonald-91}.

The description of the cohomology rings of the flag manifolds
$U(n)/T^n$and Grassmann manifolds $G_{n,k}=U(n)/(U(k)\times U(n-k))$
in the terms of Schubert polynomials is given in~\cite{F}.

The description of the {\it complex cobordism ring} of the flag manifolds
$G/T$, for $G$ compact, connected Lie group and $T$ its maximal torus,
 in the terms of the Schubert polynomials calculus is given
in~\cite{Bressler-Evens-90, Bressler-Evens-92}.
\end{rem}

\subsection{Grassmann manifolds.} As a next application we will compute
cobordism class, characteristic numbers $s_{\omega}$ and,
consequently, Chern numbers for invariant complex structure on
Grassmannian $G_{4,2} = U(4)/(U(2)\times U(2)) = SU(4)/S(U(2)\times
U(2))$. Note that, it follows by~\cite{BH} that, up to equivalence,
$G_{4,2}$ has one invariant complex structure $J$. The corresponding
Lie algebra description for $G_{4,2}$ is $A_{3}/(\TT ^{1}\oplus A_1
\oplus A_1)$.

The number of the fixed points under the canonical action of $T^{3}$
on $G_{4,2}$ is, by Theorem~\ref{number}, equal to 6. Let
$x_1,x_2,x_3,x_4$ be canonical coordinates on maximal abelian
algebra for $A_3$. Then $x_1, x_2$ and $x_3, x_4$ represents
canonical coordinates for $A_1\oplus A_1$. The weights of this
action at identity point $(T_{e}(G_{4,2}), J)$are given by the
positive complementary roots $x_1-x_3, x_1-x_4, x_2-x_3, x_2-x_4$
for $A_3$ related to $A_1\oplus A_1$ that define $J$.

The Weyl group $W_{U(4)}$ is the symmetric group of permutation on
coordinates $x_1,\ldots ,x_4$ and the Weyl group $W_{U(2)\times
U(2)} = W_{U(2)}\times W_{U(2)}$ is the product of symmetric groups
on coordinates $x_1, x_2$ and $x_3, x_4$ respectively. Let
$\rw_{j}\in  W_{U(4)}/W_{U(2)\times U(2)}$. Corollary~\ref{chom}
gives that the cobordism class $[G_{4,2}]$ is the coefficient for
$t^{4}$ in polynomial
\begin{equation}\label{gr}
\begin{split}
\sum _{j=1}^{6}\rw _{j}\Big( \frac{f(t(x_1-x_3))f(t(x_1-x_4))
f(t(x_2-x_3))
f(t(x_2-x_4))}{(x_1-x_3)(x_1-x_4)(x_2-x_3)(x_2-x_4)}\Big) = \qquad \qquad\\
= \frac{1}{4} L\Big( (x_1-x_2)(x_3-x_4)f(t(x_1-x_3))
f(t(x_1-x_4)) f(t(x_2-x_3)) f(t(x_2-x_4)) \Big)\; ,
\end{split}
\end{equation}
where $f(t) = 1 + a_{1}t + a_{2}t^2 + a_{3}t^3 + a_{4}t^4$.

Expanding formula~\eqref{gr} we get that
\begin{equation}\label{CCG42}
[G_{4,2}] = 2(3a_{1}^4+ 12a_{1}^2a_{2} + 7a_{2}^2 + 2a_{1}a_{3}-10a_4) \; .
\end{equation}
The characteristic numbers $s_{\omega}$ can be read off form this formula:
\[
s_{(4,0,0,0)}=6, \quad s_{(2,1,0,0)}=24, \quad s_{(0,2,0,0)}=14,
\quad s_{(1,0,1,0)}=4, \quad s_{(0,0,0,1)} = -20 \; .
\]
The coefficients $\beta _{\omega \xi}$ from  Theorem~\ref{sc} can be
explicitly computed and for $8$-dimensional manifold give the
following relation between characteristic numbers $s_{\omega}$ and
Chern numbers:
\[
s_{(0,0,0,1)}= c_{1}^4 - 4c_{1}^2c_{2} + 2c_{2}^2 + 4c_{1}c_{3} - 4c_{4}
, \quad s_{(2,1,0,0)} = c_1c_3 - 4c_4 \; ,
\]
\[
s_{(0,2,0,0)} = c_{2}^2 - 2c_1c_3 + 2c_4, \quad
s_{(1,0,1,0)}=c_1^2c_2-c_1c_3+4c_4-2c_2^2,\quad s_{(4,0,0,0)} = c_4 \; .
\]
We deduce that the Chern numbers for $(G_{4,2}, J)$ are
\[
c_4 = 6, \quad c_1c_3 = 48, \quad  c_2^2 = 98, \quad c_1^2c_2 = 224,
\quad c_1^4 = 512 \; .
\]
The given example  generalizes as follows. Denote by  $\Delta_{p,q}
= \prod\limits_{p\leqslant i<j\leqslant q}(x_i-x_j)$, then
$\Delta_n = \Delta_{1,n}$.

\begin{thm}\label{CGR}
The cobordism class for Grassmann manifold $G_{q+l,l}$ is given as
the coefficient for $t^{lq}$ in the series in $t$
\begin{equation}
\sum _{\sigma\in S_{q+l}/S_{q}\times S_{l}}\sigma\Big(\prod
\frac{f(t(x_i-x_j))}{(x_i-x_j)}\Big) = \frac{1}{q!l!}L \Big(
\Delta_q\Delta_{q+1,q+l}\prod f(t(x_i-x_j)) \Big)\; ,
\end{equation}
where $1 \leqslant i \leqslant q, \; (q+1) \leqslant j \leqslant
(q+l)$ and $S_{q+l}$ is the symmetric group.
\end{thm}

Let us introduce the polynomials $Q_{(q+l,l)\xi}$ defined with
\[
\Delta_q\Delta_{q+1,q+l}\prod\limits_{\underset{q+1\leqslant j\leqslant q+l}{1\leqslant i\leqslant q}}f(t(x_i-x_j))=
\sum_{|\xi | \geqslant \frac{(q+l)^2-(q+l)}{2}}Q_{(q+l,l)\xi}(a_1,\ldots , a_n,\ldots)
t^{|\xi |-\frac{(q+l)^2-(q+l)}{2}}x^{\xi}\; ,
\]
where $\xi = (j_1,\ldots ,j_{q+l})$ and $|\xi |=\sum\limits_{k=1}^{q+l}j_k$.
Appealing to Theorem~\ref{CGR} we obtain the following.

\begin{cor}\label{CDEL}
The cobordism class for Grassmann manifold $G_{q+l,l}$ is given with
\begin{equation}
[G_{q+l,l}] =\frac{1}{q!l!}\sum _{\sigma\in S_{q+l}}\sg (\sigma)Q_{(q+l,l)\sigma\delta}(a_1,\ldots,a_{ql})\; ,
\end{equation}
where $\delta = (q+l-1,q+l-2,\ldots,0)$.
\end{cor}

\begin{ex}
For $q=l=2$ the calculations give that the polynomials $Q_{\sigma\delta}=Q_{(4,2)\sigma\delta}(a_1,a_2,a_3,a_4)$ are as follows
\[
Q_{(3,2,1,0)}=-Q_{(2,3,1,0)}=-Q_{(3,2,0,1)}=Q_{(2,3,0,1)}=Q_{(1,0,3,2)}=-Q_{(1,0,2,3)}=-Q_{(0,1,3,2)}=Q_{(0,1,2,3)}\]
\[=a_1^4+4a_2^2-4a_1a_3\; ,
\]
\[
Q_{(2,1,3,0)}=-Q_{(1,2,3,0)}=-Q_{(2,1,0,3)}=Q_{(1,2,0,3)}=Q_{(3,0,2,1)}=-Q_{(3,0,1,2)}=-Q_{(0,3,2,1)}=Q_{(0,3,1,2)}\]\[=a_1^4+4a_1^2a_2+a_2^2+6a_1a_3-6a_4\; ,
\]
\[
Q_{(1,3,2,0)}=-Q_{(3,1,2,0)}=Q_{(3,1,0,2)}=-Q_{(1,3,0,2)}=Q_{(2,0,1,3)}=-Q_{(2,0,3,1)}=Q_{(0,2,3,1)}=-Q_{(0,2,1,3)}\]\[=a_1^4+8a_1^2a_2+2a_2^2-4a_4\; ,
\]
Using Corollary~\ref{CDEL} we obtain the formula~\eqref{CCG42} for the cobordism class $[G_{4,2}]$.
\end{ex}

\subsection{Homogeneous space $SU(4)/S(U(1)\times U(1)\times U(2))$}\label{M10}
Following~\cite{BH} and~\cite{KT} we know that 10-dimensional space
$M^{10} = SU(4)/S(U(1)\times U(1)\times U(2))$ admits, up to
equivalence, two invariant complex structure $J_1$ and $J_2$ and one
non-integrable invariant almost complex structure $J_{3}$. We
provide here the description of cobordism classes for all of three
invariant almost complex structures. The Chern numbers for all the
invariant almost complex structures are known and they have been
completely computed in~\cite{KT} through multiplication in
cohomology. We provide also their computation
using our method.

The corresponding Lie algebra description for $M^{10}$ is
$A_{3}/(\TT ^{2}\oplus A_1)$. Let $x_1, x_2, x_3, x_4$ be canonical
coordinates on maximal Abelian subalgebra for $A_3$. Then $x_1, x_2$
represent canonical coordinates for $A_1$. The number of fixed
points under the canonical action of $T^{3}$ on $M^{10}$ is, by
Theorem~\ref{number}, equal to 12.

\subsubsection{The invariant complex structure $J_1$.}\label{M101}
The weights of the action of $T^{3}$ on $M^{10}$ at identity point
related to $J_1$ are given by the  complementary roots $x_1-x_3,
x_1-x_4, x_2-x_3, x_2-x_4, x_3-x_4$ for $A_{3}$ related to $A_1$,
(see~\cite{BH},~\cite{KT}). The cobordism class $[M^{1}, J_1]$ is,
by Corollary~\ref{chom}, given as the coefficient for $t^{5}$ in
polynomial
\[
\sum _{j=1}^{12}\rw _{j}\Big( \frac{f(t(x_1-x_3))f(t(x_1-x_4))
f(t(x_2-x_3))f(t(x_2-x_4))
f(t(x_3-x_4))}{(x_1-x_3)(x_1-x_4)(x_2-x_3)(x_2-x_4)(x_3-x_4)}\Big ) \; ,
\]
where $\rw_{j}\in  W_{U(4)}/W_{U(2)}$ and $f(t) = 1 + a_{1}t +
a_{2}t^2 + a_{3}t^3 + a_{4}t^4 + a_{5}t^5$.

Therefore we get that
\[
[M^{10}, J_{1}]= 4(3a_1^5 + 12a_1^3a_2 + 7a_1a_2^2 -
5a_1^2a_3 - 2a_2a_3 - 10a_1a_4 + 5a_5) \; .
\]
By Theorem~\ref{sc} we get the following relations between
characteristic numbers $s_\omega$ and Chern numbers  for $(M^{10},
J_1)$.
\[
s_{(0,0,0,0,1)} = 20 = c_1^5 - 5c_1^3c_2 + 5c_1^2c_3 + 5c_1c_2^2 -
5c_1c_4 - 5c_2c_3 + 5c_5 \; ,
\]
\[
s_{(1,2,0,0,0)}= 28 = c_2c_3 - 3c_1c_4 + 5c_5, \quad
s_{(2,0,1,0,0)} = -20 = c_1^2c_3 - c_1c_4 - 2c_2c_3 + 5c_5 \; ,
\]
\[
s_{(0,1,1,0,0)} = -8 = -2c_1^2c_3 + c_1c_2^2 - c_2c_3 + 5c_1c_4 - 5
c_5, \quad s_{(3,1,0,0,0)} = 48 = c_1c_4 - 5c_5 \; ,
\]
\[
s_{(1,0,0,1,0)} = -40 = c_1^3c_2 - c_1^2c_3 - 3c_1c_2^2 + c_1c_4 +
5c_2c_3 - 5c_5, \quad s_{(5,0,0,0,0)} = 12 = c_5 \; .
\]
This implies that the Chern numbers for $(M^{10}, J_1)$ are as
follows:
\[
c_5 = 12, \quad c_1c_4 = 108, \quad c_2c_3 = 292, \quad c_1^2c_3 =
612, \quad c_1c_2^2 = 1028, \quad c_1^3c_2 = 2148, \quad c_1^5 =
4500 \; .
\]

\subsubsection{The invariant complex structure $J_2$.} The weights of the
action of $T^{3}$ on $M^{10}$ at identity point related to $J_2$ are
given by the positive complementary roots $x_4-x_1, x_4-x_2,
x_4-x_3, x_1-x_3, x_2-x_3$ for $A_{3}$ related to $A_1$,
(see~\cite{BH},~\cite{KT}). The cobordism class $[M^{1}, J_2]$ is,
by Corollary~\ref{chom}, given as the coefficient for $t^{5}$ in
polynomial
\[
\sum _{j=1}^{12}\rw _{j}\Big( \frac{f(t(x_4-x_1))f(t(x_4-x_2))
f(t(x_4-x_3))f(t(x_1-x_3))
f(t(x_2-x_3))}{(x_4-x_1)(x_4-x_2)(x_4-x_3)(x_1-x_3)(x_2-x_3)}\Big) \; ,
\]
where $\rw_{j}\in  W_{U(4)}/W_{U(2)}$ and $f(t) = 1 + a_{1}t +
a_{2}t^2 + a_{3}t^3 + a_{4}t^4 + a_{5}t^5$.

Therefore we get that
\[
[M^{10}, J_{2}]= 4(3a_1^5 + 12a_1^3a_2 + 7a_1a_2^2 - 5a_1^2a_3 +
8a_2a_3 - 10a_1a_4 - 5a_5) \; .
\]

\begin{rem}
We could proceed with the description of cobordism class $[M^{10},
J_{2}]$ appealing to the description of $J_1$ from~\ref{M101}
and applying the results from Section~\ref{eqsthom}. The corresponding description of the weights for the
action on $T^3$ on $(M^{10}, J_2)$ looks like
\[
(+1)\cdot (x_1-x_3), (-1)\cdot (x_1-x_4), (+1)\cdot (x_2-x_3),
(-1)\cdot (x_2-x_4), (-1)\cdot (x_3-x_4) \; ,
\]what means that
\[
a_{1}(e)=+1, a_{2}(e)=-1, a_{3}(e)=+1, a_{4}(e)=-1,a_{5}(e)=-1\; .
\]
Since $J_1$ and $J_2$ define on $M^{10}$ an opposite orientation, we
have that $\epsilon =-1$ and  it follows that the fixed points have
sign $+1$.
\end{rem}

Applying the same procedure as above we get that the Chern
numbers for $(M^{10}, J_{2})$ are:
\[
c_5 = 12, \quad c_1c_4 = 108, \quad c_2c_3 = 292, \quad c_1^2c_3 = 612,
\quad c_1c_2^2 = 1068, \quad c_1^3c_2 = 2268, \quad c_1^5 = 4860 \; .
\]

\subsubsection{The invariant almost complex structure $J_3$.} The weights for
the action of $T^{3}$ on $M^{10}$ at identity point related to $J_3$
are given by complementary roots $x_1-x_3, x_2-x_3, x_4-x_1,
x_4-x_2, x_3-x_4$, (see~\cite{KT}). Using Corollary~\ref{chom} we get that the
cobordism class for $(M^{10}, J_3)$ is
\[
[M^{10}, J_{3}] = 4(3a_1^5 - 12a_1^3a_2 + 7a_1a_2^2 + 15a_1^2a_3 -
12a_2a_3 -10a_1a_4 + 15a_5) \; .
\]

\begin{rem}
We could also as in the previous case proceed with the computation
of the cobordism class $[M^{10}, J_{3}]$ appealing on the on the
description of $J_1$ from~\ref{M101}. The corresponding description
of the weights is
\[
(+1)\cdot (x_1-x_3), (-1)\cdot (x_1-x_4), (+1)\cdot (x_2-x_3),
(-1)\cdot (x_2-x_4), (+1)\cdot (x_3-x_4)\; .
\]
Since $J_1$ and $J_2$ define the same orientation it follows that
all fixed points for an action $T^3$ on $M^{10}$ have sign $+1$.
\end{rem}

The characteristic numbers for $(M^{10}, J_{3})$ are given as
coefficients in its cobordism class, what, as above,  together with
Theorem~\ref{sc} gives that the Chern numbers for $(M^{10}, J_{3})$
are as follows:
\[
c_5 = 12, \quad c_1c_4 = 12, \quad c_2c_3 = 4, \quad c_1^2c_3 = 20,
\quad c_1c_2^2 = -4, \quad c_1^3c_2 = -4, \quad c_1^5 = -20 \; .
\]

\begin{rem}
Further work on the studying of Chern numbers
and the geometry for the generalizations of this example is done
in~\cite{H05} and in~\cite{KT}.
\end{rem}

\subsection{Sphere $S^{6}$.}According to~\cite{BH} we know that the
sphere $S^{6}=G_2/SU(3)$ admits $G_2$-invariant almost complex
structure, but it does not admit $G_2$-invariant complex structure.
The existence of an invariant almost complex structure follows from
the fact that $SU(3)$ is connected centralizer of an element of
order $3$ in $G_2$ which generates it's center, while the
non-existence of an invariant complex structure is because the
second Betti number for $S^6$ is zero. Note that being isotropy
irreducible, $S^6=G_2/SU(3)$ has unique, up to conjugation,
invariant almost complex structure  $J$.

The roots for the Lie algebra $\gg _{2}$ are given with (see~\cite{Onishchik})
\[
\pm x_1,\; \pm x_2,\; \pm x_3,\;  \pm (x_1-x_2),\; \pm(-x_1+x_3),\; \pm(-x_2+x_3),
\]
where $x_1+x_2+x_3=0$. It follows that the system of complementary roots for $\gg _{2}$ related to $A_2$ is
$\pm x_1$, $\pm x_2$, $\pm x_3$. According to~\ref{ac}, since $S^6=G_2/SU(3)$ is isotropy irreducible, this implies that the roots of   an existing invariant almost complex
structure $J$ on $S^{6}$ are
\[
\alpha _{1} = x_1,\; \alpha_{2}=x_2,\;
\alpha_{3}= x_3=-(x_1+x_2)\; .
\]
The canonical action of a common maximal torus $T^2$ on
$S^6=G_2/SU(3)$ has $\chi (S^6)=2$ fixed points.

By Theorem~\ref{weights} we get that the weights at
the fixed points for this action are given by the action of the Weil group $W_{G_2}$ up to the action of the Weil group $W_{SU(3)}$ on the roots for $J$:
\[
x_{1},\; x_{2},\; x_3\;\; \mbox{and}\;\; -x_{1},\;
-x_{2},\; -x_3\; .
\]
Since the weights at these two fixed points are of opposite signs,  Corollary~\ref{ch-tor-hom} implies that in the Chern character of the universal toric genus
for $(S^6, J)$ the coefficients for $a^{\omega}$ are going to be zero for  $\|\omega \|$ being even.

As the  almost complex structure $J$ is invariant under the action of the group $G_2$, it follows by Lemma~\ref{L1} that the universal toric genus for $(S^{6},J)$ belongs to the image of the ring  $U^{*}(BG_2)$ in $U^{*}(BT^2)$ of the map induced by embedding $T^2\subset G_2$ . Furthermore,
using  Corrollary~\ref{ch-tor-hom} we obtain that the universal toric genus for $(S^6,J)$ is series in $\sigma_{2}$ and $\sigma_{3}$, where $\sigma_2$ and $\sigma_3$ are the elementary symmetric functions in three variables.  The direct computation  gives that the beginning terms in the series of the Chern chraracter are
\begin{align*}
ch_{U}\Phi(S^{6},J) = & 2(a_1^3 - 3a_1a_2 + 3a_3)+\\
                      & + 2( a_1a_2^2-2a_1^2a_3 - a_2a_3+5a_1a_4-5a_5)\sigma _{2}+\\
                      & + 2(a_1a_3^2-2a_1a_2a_4-a_3a_4+2a_1^2a_5+3a_2a_5-7a_1a_6+7a_7)\sigma_{2}^{2} +\\
                      & + 2(3a_9-3a_1a_8-3a_2a_7+6a_3a_6-3a_4a_5+3a_1^2a_7-3a_1a_2a_6-3a_1a_3a_5+\\
                      & + 3a_1a_4^2+3a_2^2a_5-3a_2a_3a_4+a_3^3)\sigma_{3}^{2}+\\
                      & + 2(-9a_9+9a_1a_8-5a_2a_7+3a_3a_6-a_4a_5-2a_1^2a_7+\\
                      & + 2a_1a_2a_6-2a_1a_3a_5+a_1a_4^2)\sigma_{2}^{3}+ \ldots
\end{align*}
In particular, we obtain that the cobordism class for $(S^6,J)$ is
\begin{equation}\label{cobsph}
[S^{6},J]=2(a_1^3 - 3a_1a_2 + 3a_3) \; .
\end{equation}

\begin{rem}
We can also compute cobordism class $[S^{6}, J]$ using relations
between Chern numbers and  characteristic numbers for an invariant
almost complex structure given by Theorem~\ref{sc}:
\[
c_3 = s_{(3,0,0)},\; c_1c_2 - 3c_3 = s_{(1,1,0)},\; c_1^3 - 3c_1c_2
+ 3c_3 = s_{(0,0,1)}\; .
\]
Since for $S^6$ we obviously have that $c_1c_2=c_1^3=0$ and
$c_{3}=2$, it implies $s_{(3,0,0)}=2$, $s_{(1,1,0)}=-s_{(0,0,1)}=6$
what gives formula~\eqref{cobsph}.
\end{rem}

\subsubsection{On stable complex structures.}
If now $c_{\tau}$ is an arbitrary stable complex structure on $S^6$,
equivariant under the given action of $T^2$, than 
Corollary~\ref{stableweights} implies that the weights at the fixed points of this
action related to $c_{\tau}$ are given with
\[
a_{1}(1)\alpha_{1},\; a_{2}(1)\alpha_{2},\;a_{3}(1)\alpha_{3},\;
\mbox{and}\;
a_{1}(2)(-\alpha_{1}),\;a_{2}(2)(-\alpha_{3}),\;a_{3}(1)(-\alpha_{2})\;
.
\]
By Corollary~\ref{stablesigns} the signs at fixed points are
\[
\sg(i)=\epsilon \cdot \prod_{j=1}^{3}a_{j}(i),\;\; i=1,2\; .
\]
Corollary~\ref{stablecoeffzero} implies that the coefficients $a_{j}(i)$ should satisfy the following
equations:
\[
a_{1}(1)+a_{1}(2)=a_{2}(1)+a_{2}(2) = a_{3}(1)+a_{3}(2),\;
a_{1}(1)a_{2}(1)-a_{1}(2)(a_{1}(1)-a_{2}(1))=1 \; .
\]
These equations have ten solutions which we write as the couples of
triples that correspond to the fixed points. Two of them are
$(1,1,1), (1,1,1)$ and $(-1,-1,-1), (-1,-1,-1)$ and correspond to an
invariant almost complex structure $J$ and it's conjugate. The other
eight couples are of the form $(i,j,k), (-i,-j,-k)$ where $i,j,k=\pm
1$ and they  describe the weighs of any other stable complex
structure on $S^6$ equivariant under the given torus action. Note
that, since $\tau (S^{6})$ is trivial bundle, we have on $S^{6}$
many stable complex structures different from $J$.
The fact that, for any $T^2$-equivariant stable complex structure on $S^6$ different from $J$ or it's conjugate, the weights at two fixed points differ by sign, together with   Proposition~\ref{cobc} proves the following Proposition.

\begin{prop}
The cobordism class for $S^6$
related to any $T^2$-equivariant stable complex  structure that is not equivalent to $G_2$-invariant almost complex structure is  trivial. It particular, besides of described $G_2$-invariant almost
complex structure, $S^6$ does not admit any other almost complex
structure invariant under the  canonical action of $T^2$.
\end{prop}

\bibliographystyle{amsplain}

\end{document}